\documentclass[11pt,authoryear]{elsarticle}
\usepackage[a4paper,hmargin=1.8cm,vmargin=1.8cm]{geometry}
\usepackage{capt-of}
\usepackage{setspace}
\AtBeginDocument{
  \setlength{\abovedisplayskip}{4pt plus 2pt minus 1pt}
  \setlength{\belowdisplayskip}{4pt plus 2pt minus 1pt}
  \setlength{\abovedisplayshortskip}{3pt plus 1pt}
  \setlength{\belowdisplayshortskip}{4pt plus 1pt minus 1pt}
}
\makeatletter
\let\ps@pprintTitle\ps@plain
\makeatother
\usepackage{amsthm}
\RequirePackage{tgtermes}
\RequirePackage{newtxtext}
\RequirePackage{bm}

\usepackage{amsmath,amssymb,amsfonts}
\RequirePackage{newtxmath}
\allowdisplaybreaks[4]
\usepackage{enumerate}
\usepackage{graphicx}
\usepackage{xcolor}
\usepackage{soul}
\definecolor{ReviewYellow}{RGB}{255,245,140}

\usepackage{algorithm}
\usepackage{algpseudocode}
\usepackage{textcomp}

\biboptions{authoryear,round,sort}
\AtBeginDocument{\let\cite\citep}
\usepackage{hyperref}
\hypersetup{hidelinks}
\let\mathscr\mathcal
\makeatletter
\newenvironment{breakablealgorithm}
{
	\begin{center}
		\small
		\refstepcounter{algorithm}
		\hrule height.8pt depth0pt \kern2pt
		\renewcommand{\caption}[2][\relax]{%
			{\raggedright\textbf{\ALG@name~\thealgorithm} ##2\par}%
			\ifx\relax##1\relax
			\addcontentsline{loa}{algorithm}{\protect\numberline{\thealgorithm}##2}%
			\else
			\addcontentsline{loa}{algorithm}{\protect\numberline{\thealgorithm}##1}%
			\fi
			\kern2pt\hrule\kern2pt
		}
	}{
		\kern2pt\hrule\relax
	\end{center}
}
\makeatother

\theoremstyle{plain}
\newtheorem{theorem}{Theorem}[section]
\newtheorem{lemma}[theorem]{Lemma}
\newtheorem{proposition}[theorem]{Proposition}
\newtheorem{corollary}[theorem]{Corollary}
\newtheorem{assumption}{Assumption}[section]
\newtheorem{condition}[theorem]{Condition}
\theoremstyle{definition}
\newtheorem{definition}[theorem]{Definition}
\newtheorem{remark}[theorem]{Remark}

\newcommand{\cP}{\mathcal{P}}

\begin{document}
\begin{frontmatter}
\title{\resizebox{\textwidth}{!}{Approximation and Error Bounds for Zero-Sum Stopping Games with Probability Criteria}}
\author[sysugz]{Xin Wen}
\ead{wenx78@mail.sysu.edu.cn}
\author[sysugz]{Xianping Guo}
\ead{mcsgxp@mail.sysu.edu.cn}
\author[sysusz]{Xin Guo}
\ead{guox87@mail.sysu.edu.cn}
\author[cuhk]{Zewu Zheng\corref{cor1}}
\ead{zhengzw@link.cuhk.edu.hk}
\cortext[cor1]{Corresponding author.}
\address[sysugz]{School of Mathematics, Sun Yat-sen University, Guangzhou 510275, China}
\address[sysusz]{School of Science, Sun Yat-sen University, Shenzhen, China}
\address[cuhk]{Department of Statistics, The Chinese University of Hong Kong, China}
\begin{abstract}
This paper studies a two-person zero-sum stochastic game, which must terminate
when either player quits the game or a state of the game first enters a prescribed target set. The performance
criterion is the probability that the total reward reaches a given threshold. First, under
a mild condition, we prove the existence of the value and a Nash equilibrium of the
game by developing a comparison theorem. Second, for solving the game with infinite
states, we develop a finite-state approximation method based on truncation, and derive
explicit error bounds between the values of the truncated and original games. Third,
using the error bounds, we design an effective algorithm, by which we can not only
obtain any approximation of the value and Nash equilibria of the original game, but
also provide convergence analysis of the approximation. Finally, a queueing example is
used to illustrate the advantage of our condition and the effectiveness of our algorithm.
\end{abstract}
\end{frontmatter}

\section{Introduction}
\label{sec:introduction}

Stochastic games provide a powerful framework for modeling multi-agent decision-making in stochastic dynamic environments \cite{Shapley1953Stochastic, Neyman2003,JaskiewiczNowak2018}. In the classical settings, the agents typically optimize expected cumulative payoffs under the risk-neutral or risk-sensitive criteria. However, in many risk-sensitive applications, the agents may follow  the consideration that security takes priority over benefit. Based on this consideration, the agents are not concerned with the maximization of their expected total rewards, instead, they hope to maximize the likelihood of their rewards reaching  a given threshold, in which the likelihood is denoted by probability, for example, the probability that the rewards of the agents reaches any given target level, the probability that the surplus of an insurance company reaches a given threshold before  ruin occurs, and so on. These situations have motivated the study of probability criteria in stochastic games, where the objective is to maximize the probability that the total reward reaches a given threshold, or equivalently to minimize the probability that it falls below that threshold. Compared with the risk-neutral and risk-sensitive expected criteria for cumulative payoffs \cite{WeiChen2019,CavazosHernandez2019,ChenGuo2023}, the probability criterion directly captures risk performance and has practical importance.

For stochastic games under the probability criterion, several works have been developed. Precisely, \citet{HGP-2017} established the Shapley equation for discrete-time zero-sum games under the first-passage probability criterion, and proved the existence of the value and  a pair of optimal policies. \citet{HGW-2023} studied finite-horizon probability criterion for semi-Markov games, and developed a value iteration algorithm, while \citet{Bhabak2022} and \citet{Bhabak2021} extended the results of \citet{HGP-2017} to semi-Markov and continuous-time settings, respectively. \citet{Wei2025} extended the fully observable discrete-time model of \citet{HGP-2017} to partially observable Markov games. In all these studies, the transition probabilities and reward functions are independent of time. By contrast, \citet{GuoWen2025Nonstationary} studied nonstationary nonzero-sum Markov games in which both may depend on time. It should be noted that all these works \cite{HGP-2017,HGW-2023,Bhabak2022, Bhabak2021,Wei2025,GuoWen2025Nonstationary} assume that the game cannot stop and the  agents  can only influence the game evolution through their policies and do not have the option to quit the game.

In practice, a game may actively terminate when an agent of the  game  quits the game to avoid adverse outcomes, or a game must passively stop when a state of the game  enters some critical region, such as the  bankruptcy or equipment failure state. In all these scenarios, the stopping decision directly affects the probability of reaching a target, making the probability criterion particularly relevant. To the best of our knowledge, such games when the active stopping and passive termination coexist, have not yet been systematically studied under the probability criterion.
For the probability criterion, a closely related work is given by \citet{Zhang2026}, who studied a fixed finite-horizon and special stopping game, in which only agents of the games are allowed to quit actively, but
the game does not terminate passively. Therefore, it is of both theoretical and practical importance to establish the existence of the value and Nash equilibria in stopping games under the probability criterion, where active stopping and passive termination coexist. This is the first motivation of the study in this paper.

A further challenge in stochastic games is computational. For the games with infinite states, direct computation of the values and Nash equilibria is generally infeasible. The existing approximating algorithms
\cite{HGW-2023,GuoWen2025Nonstationary,Wen2022Optimal} for the Markov games and decision processes with passive
termination under the probability criterion generally require a finite-state representation
for direct numerical implementation, and those in \citet{GuoWen2025Nonstationary,Wen2022Optimal} rely on stronger conditions; see Remark~\ref{rem:weaker_condition}
for details. Moreover, \citet{HGW-2023,GuoWen2025Nonstationary} do not provide a convergence analysis of
approximate Nash equilibria, and \citet{Zhang2026} did not develop a numerical algorithm. To handle computation on countably infinite
state spaces, a natural approach is to approximate the original game by
finite-state games. Related finite-approximation results have been established
in \citet{DufourPrietoRumeau2019,TidballAltman1996} for risk-neutral zero-sum discounted games, and in
\citet{LorenzoHernandezPrietoRumeau2015,TidballAltman1996,TidballPourtallierAltman1997}
for risk-neutral zero-sum average games. The computation methods for games with risk-neutral criteria
require further analysis before being applied to the present game,
which has no discount factor and does not assume a unichain condition.
Moreover, error estimates are needed for the approximations
considered here.
To give an effective algorithm for computing the values and Nash equilibria  with the probability criterion is another motivation for this paper.

In this paper, we study the  stopping games under  the probability criterion, and aim to establish the existence and computation of the value and Nash equilibria  for the games. Our contributions are as follows:

First, under a general condition, we establish a comparison theorem for the probability criterion (Theorem~\ref{thm:relationship} below), which forms the technical basis for the subsequent existence and computation results. To verify the general condition, we further introduce a Lyapunov condition, which is not only weaker than the conditions used in \citet{HG-2020,GuoWen2025Nonstationary,Wen2022Optimal,Huang2011,HGP-2017,Liu2001}, but also is sufficient for the general condition above (see Remark~\ref{rem:weaker_condition}). Moreover, we derived a Shapley operator for the game with stopping and the probability criterion, by which we construct two monotone sequences of functions,  and prove that they converge to the unique solution to the Shapley equation. Using the comparison theorem and the Shapley equation, we prove the existence of the value and a Nash equilibrium (Theorem~\ref{game-Nash}), and also show that the  difference of the two sequences above can be used to provide an online stopping rule of computing the  value and a Nash equilibrium of the games.

Second, to handle computation of the values and Nash equilibria for the stopping games with infinite states, we introduce a truncation approximation method. Under the Lyapunov condition, we construct a sequence of truncated games, verify that every truncated game inherits this Lyapunov condition, and therefore each truncated game has a value and a Nash equilibrium (see Corollary~\ref{cor:truncated-game-Nash} below). Moreover, using Doob's maximal inequality for nonnegative supermartingales, we derive explicit error bounds between the values of the truncated games and that of the original game (see Theorem \ref{thm:truncation} below). These error bounds are new
(see Remark~\ref{rem:truncation-novelty}), and they gradually vanish as the number of the states of truncated games increases accordingly.
Consequently, the values of the truncated games are proven to converge to the value of the original game.

Third, using the explicit error bounds obtained here, we develop an  effective value iteration algorithm
(Algorithm~\ref{alg:truncVI}), by which we can obtain any approximation of the values and a Nash equilibrium of the truncated and original games  in finitely many iterations
(see Theorems~\ref{thm:truncated-convergence} and \ref{thm:main_epsilon} below). Moreover, we have found a suitable online stopping rule given by the  difference between the two sequences above, and also provided rigorous  convergence analysis of the algorithm. Precisely, for each fixed truncated
game, as the accuracy decreases to zero, we prove that the (accuracy-dependent) approximate values obtained by Algorithm~\ref{alg:truncVI} converge to the value of  the truncated game, and also show that any accumulation point of the approximate equilibria of the truncated game, forms a Nash equilibrium of the truncated game
(see Theorem~\ref{thm:truncated-convergence} below). Furthermore, as the accuracy decreases to zero and the number of states of the truncated games increases
accordingly, we also prove that   every accumulation point of accuracy-dependent approximate Nash equilibria of the truncated games is a Nash equilibrium  of the original game (see Theorem~\ref{thm:main_epsilon} below).

Finally, to illustrate the advantage
of our condition and the effectiveness of our algorithm, we present a controlled queueing example, for which our Lyapunov condition is satisfied, while some of the conditions in \citet{HG-2020,GuoWen2025Nonstationary,Wen2022Optimal,Huang2011,HGP-2017,Liu2001} fail to hold, and for which some numerical results are obtained by Algorithm~\ref{alg:truncVI}.

The paper is organized as follows.
 Section~2 describes the game model and the probability criterion.
 Section~3 gives suitable assumptions and the comparison theorem.
 Section~4 establishes the existence of the value and a Nash equilibrium.
 Section~5 develops the truncation approximation and establishes explicit error bounds.
 Section~6 presents the value-iteration algorithm and its convergence analysis of both the truncated and the original games.
 Section~7 provides the numerical example.

\section{Game Model}\label{sec:model}
	
	We first introduce some notations used throughout this paper.
	Given a Borel space $E$, denote by $\mathcal{B}(E)$ the Borel $\sigma$-algebra of $E$, by $\mathcal{P}(E)$ the family of probability measures on $\mathcal{B}(E)$, by $\delta_z$ the Dirac measure at $z\in E$, and by $\mathbf{1}_C$ the indicator function of a set $C$.
In addition, $\mathbb{R}_+:=\allowbreak [0,\allowbreak \infty)$, $\mathbb{N}:=\allowbreak \{0,\allowbreak 1,\allowbreak 2,\allowbreak \ldots\}$, and $\inf\emptyset:=\allowbreak \infty$. For any $a,\allowbreak b\in\mathbb{R}$, write $a\vee b:=\allowbreak \max\{a,\allowbreak b\}$, $a\wedge b:=\allowbreak \min\{a,\allowbreak b\}$, and $a^+:=\allowbreak a\vee0$.
	Moreover, let
 $\sum_{t=\allowbreak l}^m y_t :=\allowbreak  0$
 for any sequence $\{y_t\}$ with $ m < l $.

	A two-person zero-sum first-passage stochastic game with stopping
is specified by
\begin{eqnarray}\label{def:game}
\mathcal G:=\{S,D,A,B,(A(i),B(i))_{i\in S},p,r,g^1,g^2,\lambda^*\},
\end{eqnarray}
where $S$ is a denumerable state space endowed with its discrete
$\sigma$-algebra, and $D\subset S$ is a prescribed target set
with complement $D^c:=S\setminus D$.
The action spaces of player~1 (the maximizer) and player~2
(the minimizer) are $A$ and $B$, respectively; for each $i\in S$,
the admissible action sets $A(i)\subseteq A$ and $B(i)\subseteq B$
are nonempty and finite.
The stochastic kernel $p(\cdot\mid i,a,b)$ on $S$ given
$\mathbb K:=\{(i,a,b):i\in S,\ a\in A(i),\ b\in B(i)\}$
specifies the transition probabilities.
The one-stage reward $r:\mathbb K\to\mathbb R_+$ for player~1
is the cost for player~2 and satisfies $r(i,a,b)>0$ whenever
$i\in D^c$.
The functions $g^1,g^2:S\to\mathbb R_+$ are terminal rewards
for player~1 (terminal costs for player~2), and $\lambda^*>0$
is a given initial threshold.
	
	 We now describe the evolution of the discrete time stochastic game with stopping. At each time $t=\allowbreak 0,\allowbreak 1,\allowbreak \dots,\allowbreak $ suppose that the system is in state $i_t\in D^c$ and maintains a remaining threshold $\lambda_t$, with {$\lambda_0:=\allowbreak \lambda^*$. Based on the current state $i_t\in D^c$ and the current threshold  $\lambda_t$, the following happens:
		 \begin{enumerate}
	 	\item[(i)] If player~1 decides to  quit the game or both players quit at the same time, then player 1 receives a reward $g^1(i_t)$ from player 2 and the game ends.
	 	\item[(ii)] If player~2 decides to quit the game while player~1 continues, then player~1 receives reward $g^2(i_t)$ from player 2 and the game ends.
	 	\item[(iii)] If neither player quits, the two players choose their actions $a_t\in A(i_t)$ and $b_t\in B(i_t)$ simultaneously and independently. Player 1 then receives a reward $r(i_t,\allowbreak a_t,\allowbreak b_t)$ (i.e., a cost of player 2), and the system jumps to state $i_{t+1}$ according to the transition probability $p(i_{t+1}|i_t,\allowbreak a_t,\allowbreak b_t)$. The threshold is updated as $\lambda_{t+1} :=\allowbreak \lambda_t-r(i_t,\allowbreak a_t,\allowbreak b_t)$.
	 	\item[(iv)] If the state $i_{t+1}$ is in the target set (i.e., $i_{t+1}\in D$), then the game terminates automatically; otherwise, the game continues from $(i_{t+1},\allowbreak \lambda_{t+1})$.
	 \end{enumerate}
	 Thus, for each $t\geq 0$, we can define a history up to time $t$ as
	 $h_t :=\allowbreak  (i_0,\allowbreak \lambda_0,\allowbreak a_0,\allowbreak b_0,\allowbreak \dots,\allowbreak i_{t-1},\allowbreak \lambda_{t-1},\allowbreak a_{t-1},\allowbreak b_{t-1},\allowbreak i_t,\allowbreak \lambda_t)$,
	 with $h_0=\allowbreak (i_0,\allowbreak \lambda_0)$.
Let $H_t$ be the set of finite histories up to time $t$,
endowed with the $\sigma$-algebra, and $\Lambda^*$ be the reachable threshold set, that is,  $\Lambda^*:=\allowbreak
	 \{\lambda^{*}-\sum_{k=\allowbreak 0}^{m-1}r(i_k,\allowbreak a_k,\allowbreak b_k): m\ge 0,\allowbreak  i_k\in D^c,\allowbreak  (i_k,\allowbreak a_k,\allowbreak b_k)\in\mathbb{K},\allowbreak  k=\allowbreak 0,\allowbreak 1,\allowbreak \ldots,\allowbreak m-1\},\allowbreak $ which is denumerable.

	 To regulate the action selection and formulate the optimization criterion, we introduce stopping times and policies for players.
	 First, let $
\Omega
:=\allowbreak H_{\infty}$,
and $\mathcal F$ be the corresponding $\sigma$-algebra.  For each $t\in\mathbb{N}$, the projections $i_t$, $\lambda_t$, $a_t$ and $b_t$ from $\Omega$ to the sets $S$, $\Lambda^*$, $A$ and $B$ are defined as follows.  For each path $\omega:=\allowbreak  (i_0,\allowbreak \lambda_0,\allowbreak a_0,\allowbreak b_0,\allowbreak \ldots,\allowbreak i_t,\allowbreak \lambda_t,\allowbreak a_t,\allowbreak b_t,\allowbreak \ldots)\in\Omega$,
	 \begin{eqnarray}\label{def:rv}
	 i_t(\omega):=i_t,\lambda_t(\omega):=\lambda_t, a_t(\omega):=a_t, b_t(\omega):=b_t,
	  z_0:=(i_0,\lambda_0), z_{t+1}:=(z_t,a_t,b_t,i_{t+1},\lambda_{t+1}).
	 \end{eqnarray}
Here, when no confusion arises, we use the same notation for the coordinate random variables $i_t,\lambda_t,a_t,b_t$ and their values (realization).
	
\begin{definition}\label{def:stoppingtime}
	 	A random variable $\tau:\Omega\to\mathbb{N}\cup\{\infty\}$ is called a stopping time if
	 	$\{\tau=\allowbreak t\}\in\sigma(z_t)$ for all $t\in\mathbb{N}$, where $\sigma(z_t)$ denotes the $\sigma$-algebra generated by $z_t$.
	 	The set of all stopping times is denoted by $\Gamma$.
	 \end{definition}	
	 \begin{definition}\label{def:control}
	 	\begin{enumerate}
	 		\item[(a)] A policy $\pi^1:=\allowbreak \{\pi_t^1,\allowbreak t\geq 0\}$ for player 1 is a sequence of stochastic kernels $\pi_t^1$ on $A$ given $H_t$ such that $\pi_t^1(A(i_t)|h_t)=\allowbreak 1$ for each $h_t\in H_t$, $t\in\mathbb{N}$. The family of all policies for player 1 is denoted by $\Pi^1$.
	 		\item[(b)]  A policy $\pi^1=\allowbreak \{\pi_t^1,\allowbreak t\geq0\}$ for player 1 is called stationary if there is a stochastic kernel $\phi$ on $A$ given $S\times\Lambda^*$  such that
	 		$\pi_t^1(\cdot|h_t)=\allowbreak \phi(\cdot|i_t,\allowbreak \lambda_t)$ for each $h_t\in H_t$ and $t\in\mathbb{N}$.
	 		In this case, we write $\pi^1=\allowbreak \{\pi_t^1,\allowbreak t\geq0\}$ as $\phi$ for simplicity. Denote the set of all stationary policies for player 1 by $\Pi^1_{\mathrm{s}}$. Obviously, $\Pi^1_{\mathrm{s}}\subset\Pi^1$.
	 		\item[(c)] A strategy for player 1 is a pair $(\pi^1,\allowbreak \tau^1)\in\Pi^1\times\Gamma$.
	 	\end{enumerate}
	 	Similarly, we can define a policy and a stationary policy for player 2 by replacing $A$ with $B$  and $A(i)$ with $B(i)$ in Definition \ref{def:control} (a) and (b). Denote the sets of all policies and all stationary policies for player 2 by $\Pi^2$ and $\Pi^2_{\mathrm{s}}$, respectively. Then, a strategy for player 2 is a pair $(\pi^2,\allowbreak \tau^2)\in\Pi^2\times\Gamma$.
	 \end{definition}	
	
	 For any initial state distribution $\nu\in\mathcal{P}(S)$,
initial threshold $\lambda\in\Lambda^*$, and policies $\pi^1=\allowbreak \{\pi_t^1,\allowbreak t\geq 0\}\in\Pi^1$, $\pi^2=\allowbreak \{\pi_t^2,\allowbreak t\geq 0\}\in\Pi^2$, the Ionescu-Tulcea theorem \cite{Hernandez1999Further} guarantees a unique probability measure
$\mathbb{P}_{\nu\otimes\delta_\lambda}^{\pi^1,\allowbreak \pi^2}$ on $(\Omega,\allowbreak \mathcal{F})$ satisfying:
\begin{eqnarray}\label{eq:transition-law}
	\begin{aligned}
		&\mathbb{P}_{\nu\otimes\delta_\lambda}^{\pi^1,\pi^2}
		(i_0(\omega)=i,\lambda_0(\omega)=\lambda')
		=\nu(i)\mathbf{1}_{\{\lambda'=\lambda\}},~\qquad\qquad
		\mathbb{P}_{\nu\otimes\delta_\lambda}^{\pi^1,\pi^2}
		(a_t(\omega)=a,b_t(\omega)=b|z_t)
		=\pi_t^1(a|z_t)\pi_t^2(b|z_t),\\
		&\mathbb{P}_{\nu\otimes\delta_\lambda}^{\pi^1,\pi^2}
		(i_{t+1}(\omega)=j| z_t,a_t,b_t)
		=p(j|i_t,a_t,b_t),\qquad
		\mathbb{P}_{\nu\otimes\delta_\lambda}^{\pi^1,\pi^2}
		\bigl(\lambda_{t+1}(\omega)
		=\lambda_t-r(i_t,a_t,b_t)| z_t,a_t,b_t\bigr)=1
	\end{aligned}
\end{eqnarray}
	 for  $i,\allowbreak j \in S$, $a\in A(i)$, $b\in B(i)$, and $z_t:=(i_0,\lambda_0,a_0,b_0,\cdots, i_t,\lambda_t),$ $t\geq 0$.
	 The expectation operator with respect to $\mathbb{P}_{\nu\otimes\delta_\lambda}^{\pi^1,\allowbreak \pi^2}$ is denoted by $\mathbb{E}_{\nu\otimes\delta_\lambda}^{\pi^1,\allowbreak \pi^2}$. For each $(i,\allowbreak \lambda)\in S\times\Lambda^*$, set $\mathbb{P}_{(i,\allowbreak \lambda)}^{\pi^1,\allowbreak \pi^2}:=\allowbreak \mathbb{P}_{\delta_i\otimes\delta_\lambda}^{\pi^1,\allowbreak \pi^2}$ and $\mathbb{E}_{(i,\allowbreak \lambda)}^{\pi^1,\allowbreak \pi^2}:=\allowbreak \mathbb{E}_{\delta_i\otimes\delta_\lambda}^{\pi^1,\allowbreak \pi^2}$.
	
For the given target set $D$, we denote the first-passage time by
	 $\tau_D :=\allowbreak  \inf\{t\ge0: i_t\in D\}$. Then,  for $\tau^1,\allowbreak  \tau^2\in\Gamma$,  the termination time of the game is denoted by
	\begin{equation}\label{termination time}
	\tau:=\tau_D\wedge \tau^1\wedge \tau^2.
	\end{equation}
	The total reward variable $X$ accumulated up to termination is
	\begin{eqnarray}\label{def:total-reward}
	 	X(\tau^1,\tau^2):=  \sum_{t=0}^{\tau-1} r(i_t,a_t,b_t) +
	 		\mathbf{1}_{\{\tau < \tau_D\}}\times(\mathbf{1}_{\{ \tau^1 \le \tau^2\}}g^1(i_{\tau})+\mathbf{1}_{\{ \tau^2< \tau^1\}}g^2(i_{\tau})).
	 \end{eqnarray}

	 We now define the optimization criterion. For any strategy pair $(\pi^1,\allowbreak \tau^1;\pi^2,\allowbreak \tau^2)\in(\Pi^1\times\Gamma)\times(\Pi^2\times\Gamma)$, the first-passage probability criterion is defined by
	 \begin{eqnarray}\label{def:Jnu}
	 	J_{\nu}(\pi^1,\tau^1;\pi^2,\tau^2)
	 	:=\mathbb{P}_{\nu\otimes\delta_{\lambda^*}}^{\pi^1,\pi^2}\big(X(\tau^1,\tau^2)\ge\lambda^*\big),
	 \end{eqnarray}
	which, together with $\tau_D=\allowbreak 0$ when $i_0\in D$, implies
	 \begin{eqnarray}\label{eq:Jnu-decomposition}
	 J_{\nu}(\pi^1,\tau^1;\pi^2,\tau^2)=\sum_{i\in D^c}\mathbb{P}_{(i,\lambda^*)}^{\pi^1,\pi^2}\big(X(\tau^1,\tau^2)\ge\lambda^*\big)\nu(i).
	 \end{eqnarray}
	The quantity $J_\nu(\pi^1,\allowbreak \tau^1; \pi^2,\allowbreak \tau^2)$ is the probability that the total reward reaches the prescribed threshold $\lambda^*$ with the initial state distribution $\nu$.
	The threshold $\lambda^*$ may be prescribed by a contract designer, regulator, or risk manager. Since player~2's cost is exactly player~1's reward, the event $\{X(\tau^1,\allowbreak \tau^2)\ge\lambda^*\}$ has opposite meanings for the two players: it is target attainment for player~1 and cost overrun for player~2. Consequently, player~1 aims to maximize this probability, whereas player~2 aims to minimize it.
	
The lower and upper values of the game  are defined by
	 \begin{eqnarray*}
	 	\underline{J}_{\nu}
	 	:=&\sup\limits_{(\pi^1,\tau^1)\in\Pi^1\times\Gamma}
	 	\inf\limits_{(\pi^2,\tau^2)\in\Pi^2\times\Gamma}
	 	J_{\nu}(\pi^1,\tau^1;\pi^2,\tau^2),\quad
	 	\overline{J}_{\nu}
	 	:=&\inf\limits_{(\pi^2,\tau^2)\in\Pi^2\times\Gamma}
	 	\sup\limits_{(\pi^1,\tau^1)\in\Pi^1\times\Gamma}
	 	J_{\nu}(\pi^1,\tau^1;\pi^2,\tau^2).
	 \end{eqnarray*}
	 If $\underline{J}_{\nu}=\allowbreak \overline{J}_{\nu}$, the common number is called the value of the game with respect to $(\nu,\allowbreak \lambda^*)$, denoted by $J_{\nu}^*$.
	
	 \begin{definition}\label{def:epsilon}
	 	Suppose the value $J_{\nu}^*$ exists, and fix $\epsilon\geq0$.
	 	\begin{enumerate}
	 		\item[(a)] A strategy $(\pi_{\epsilon}^1,\allowbreak \tau_{\epsilon}^1)\in\Pi^1\times\Gamma$ is called $\epsilon$-optimal for player 1 if for all $(\pi^2,\allowbreak \tau^2)\in\Pi^2\times\Gamma$,
	 		\[
	 		J_{\nu}(\pi_{\epsilon}^1,\tau_{\epsilon}^1;\pi^2,\tau^2)
	 		\geq J_{\nu}^*-\epsilon.
	 		\]
	 		\item[(b)] A strategy $(\pi_{\epsilon}^2,\allowbreak \tau_{\epsilon}^2)\in\Pi^2\times\Gamma$ is called $\epsilon$-optimal for player 2 if for all $(\pi^1,\allowbreak \tau^1)\in\Pi^1\times\Gamma$,
	 		\[
	 		J_{\nu}(\pi^1,\tau^1;\pi_{\epsilon}^2,\tau_{\epsilon}^2)
	 		\leq J_{\nu}^*+\epsilon.
	 		\]
	 		\item[(c)] If $(\pi_{\epsilon}^1,\allowbreak \tau_{\epsilon}^1)$ and $(\pi_{\epsilon}^2,\allowbreak \tau_{\epsilon}^2)$ are $\epsilon$-optimal for player 1 and player 2, respectively,
			then $(\pi_{\epsilon}^1,\allowbreak \tau_{\epsilon}^1;\pi_{\epsilon}^2,\allowbreak \tau_{\epsilon}^2)$ is called an $\epsilon$-Nash equilibrium with respect to $(\nu,\allowbreak \lambda^*)$.
	 		\item[(d)] A number $J_{\nu,\allowbreak \epsilon}$ is an $\epsilon$-approximation value of the game if $|J_{\nu,\allowbreak \epsilon}-J_{\nu}^*|\leq\epsilon$.
	 	\end{enumerate}
	 \end{definition}
	 This paper aims to establish the existence of the value and a Nash equilibrium of the game with respect to a prescribed pair $(\nu,\allowbreak \lambda^*)$, and to develop an algorithm for computing  $\epsilon$-approximation values and $\epsilon$-Nash equilibria of the game.
	 \section{Assumption and Comparison Theorem}\label{sec:assumption}
	   In this section we establish a comparison theorem that forms the technical basis for the existence of the value and Nash equilibria. Our analysis starts from the following general condition, commonly used in  \cite{GuoWen2025Nonstationary,HGP-2017,HG-2020,Wen2022Optimal,Huang2011,Liu2001}.
	 \begin{assumption}\label{ass:absorb}
	 	$\mathbb{P}_{(i,\allowbreak \lambda)}^{\pi^1,\allowbreak \pi^2}(\tau_D<\infty)=\allowbreak 1$ for every $(i,\allowbreak \lambda)\in D^c\times  \Lambda^*$ and $(\pi^1,\allowbreak \pi^2)\in \Pi^1\times \Pi^2$.
	 \end{assumption}

 To verify Assumption \ref{ass:absorb}, inspired by the
conditions in \citet{GuoSongZhang2014,PrietoRumeauHernandezLerma2012}, we give a sufficient condition weaker than those in \citet{HG-2020,HGP-2017,Huang2009,Huang2011,Liu2001}.
\begin{proposition}[Lyapunov condition]\label{prop:attenuation}
	Suppose there exist a Lyapunov function $w:S\to[1,\allowbreak \infty)$  and a constant $\rho \in [0,\allowbreak 1)$ such that
	\begin{eqnarray}\label{eq:drift_section}
	\sum_{j\in D^c} p(j|i,a,b)w(j)\leq\rho w(i)  \ \ \ \forall i\in D^c, a \in A(i), b \in B(i).
	\end{eqnarray}
	Then Assumption~\ref{ass:absorb} holds.
\end{proposition}
\begin{proof}
	For any $t\ge 0$, if $\tau_D > t$, we have $i_t \in D^c$, which together with \eqref{eq:drift_section} and \eqref{eq:transition-law}, implies that
	for every $(i,\allowbreak \lambda)\in D^c\times  \Lambda^*$ and $(\pi^1,\allowbreak \pi^2)\in \Pi^1\times \Pi^2$,
	\begin{eqnarray}\label{eq:drift_section-1}
\rho w(i_t)\mathbf{1}_{\{\tau_D > t\}}
 \ge \sum_{j \in D^c} p(j | i_t, a_t, b_t) \, w(j)\mathbf{1}_{\{\tau_D > t\}}
= \mathbb{E}_{(i,\lambda)}^{\pi^1,\pi^2}[w(i_{t+1})\mathbf{1}_{\{\tau_D > t+1\}}\mathbf{1}_{\{\tau_D > t\}} | z_t,a_t,b_t].
	\end{eqnarray}
Taking expectations gives
\begin{eqnarray}\label{prop:drift-1}
\mathbb{E}_{(i,\lambda)}^{\pi^1,\pi^2}[w(i_{t+1})\mathbf{1}_{\{\tau_D > t+1\}}]
\leq\rho\mathbb{E}_{(i,\lambda)}^{\pi^1,\pi^2}[w(i_t)\mathbf{1}_{\{\tau_D > t\}}]\leq\cdots
\leq \rho^{t+1}\mathbb{E}_{(i,\lambda)}^{\pi^1,\pi^2}[w(i_0)\mathbf{1}_{\{\tau_D > 0\}}]=\rho^{t+1}w(i),
\end{eqnarray}
which, together with $w(i) \ge 1$ for all $i \in D^c$, gives
	\begin{eqnarray}\label{prop:drift-2}
		\mathbb{P}_{(i,\lambda)}^{\pi^1,\pi^2}(\tau_D > t+1)
		=\mathbb{E}_{(i,\lambda)}^{\pi^1,\pi^2}[\mathbf{1}_{\{\tau_D > t+1\}}]\leq \rho^{t+1}w(i).
	\end{eqnarray}
	Letting $t \to \infty$ yields $	\mathbb{P}_{(i,\allowbreak \lambda)}^{\pi^1,\allowbreak \pi^2}(\tau_D =\allowbreak  \infty) =\allowbreak  0$, i.e., $	 \mathbb{P}_{(i,\allowbreak \lambda)}^{\pi^1,\allowbreak \pi^2}(\tau_D < \infty) =\allowbreak  1$.
	Thus,
	Assumption \ref{ass:absorb} holds.
\end{proof}

\begin{remark}\label{rem:weaker_condition}
\begin{itemize}

\item[(a)]
Condition~\eqref{eq:drift_section} is a Lyapunov-type condition, variants of which have been extensively used in the literature on Markov decision processes  under expected reward criteria; see, e.g.,
\citet[Assumption~2.1]{GuoSongZhang2014} and
\citet[Assumption~3.1(a)]{PrietoRumeauHernandezLerma2012}.
Here, we apply condition~\eqref{eq:drift_section} to a stopping game under a
first-passage probability criterion.
\item[(b)] It is easy to verify that the Lyapunov condition~\eqref{eq:drift_section} is weaker than the following
contraction condition:
there exists a constant $\beta>0$ such that
\begin{eqnarray}\label{one-step-absorb}
p(D\mid i,a,b)\ge\beta
\quad\text{for all }i\in D^c,\ a\in A(i),\ b\in B(i),
\end{eqnarray}
which has been used as Condition 1 in \citet{HGP-2017}, the condition for Proposition~3.1 in \citet{Huang2011}, Condition 1 in \citet{HG-2020},  and the condition for Corollary~3.2 in \citet{Liu2001}. This condition~\eqref{one-step-absorb}  implies that the operators in \citet{Huang2011,HG-2020,Liu2001} are all contractive.
\end{itemize}
\end{remark}
To achieve our aims, we next  introduce some operators. To do so,let $\mathcal{U}$ be the set of functions $u:S\times\Lambda^*\to[0,\allowbreak 1]$ satisfying
	  \begin{eqnarray}\label{def:U-boundary}
	  u(i,\lambda)=
	  \begin{cases}
1, & i\in S,\  \lambda\le 0,\\
0, & i\in D,\ \lambda>0.
\end{cases}
	  \end{eqnarray}
	  Then we define the operators
	 $\mathcal{T}$ and $\mathcal{L}$  on $\mathcal{U}$ as follows: for any $u\in\mathcal{U}$, $i\in D^c$, $\lambda\in \Lambda^*$,
	 \begin{eqnarray}\label{Tu}
\mspace{-10mu}\mathcal{T}u(i,\lambda):=\mspace{-10mu}\sup_{\phi\in\cP(A(i))}\inf_{\psi\in\cP(B(i))} \sum_{a\in A(i)}\sum_{b\in B(i)}\mathbb{G}u(i,\lambda;a,b)\phi(a)\psi(b),
	 \end{eqnarray}
	 where
\begin{eqnarray}\label{G}
	\mathbb{G}u(i,\lambda;a,b)
	:= \sum_{j\in S}u(j,\lambda-r(i,a,b))p(j|i,a,b).
\end{eqnarray}
 We further define the Shapley operator $\mathcal{L}$ as follows: for $\lambda\ge0$,
	 \begin{eqnarray}\label{Lu}
	 	\mspace{-20mu}\mathcal{L}u(i,\lambda):= \max\{ \mathbf{1}_{\{g^1(i) \ge \lambda\}},\min\{ \mathbf{1}_{\{g^2(i) \ge \lambda\}},\mathcal{T}u(i,\lambda) \} \}
	 \end{eqnarray}
	and $\mathcal{L}u(i,\allowbreak \lambda):=\allowbreak \mathcal{T}u(i,\allowbreak \lambda)=\allowbreak 1$ for $\lambda<0$. For $i\in D$, set $\mathcal{L}u(i,\allowbreak \lambda):=\allowbreak u(i,\allowbreak \lambda)$.
	
	 Since $A(i)$ and $B(i)$ are finite, $\mathcal{L}$ maps $\mathcal{U}$ into itself. For each $u\in\mathcal{U}$, we define
	 \begin{eqnarray}
	 	\tau_u^1&:=& \inf\{t \in \mathbb{N}: \mathbf{1}_{\{g^1(i_t)\ge \lambda_t\}} \geq \mathcal{L}u(i_t,\lambda_t)\}, \label{tau1}\\
	 	\tau_u^2&:=& \inf\{t \in \mathbb{N}: \mathbf{1}_{\{g^2(i_t)\ge \lambda_t\}} \leq \mathcal{L}u(i_t,\lambda_t)\}. \label{tau2}
	 \end{eqnarray}
	
	 Since the events $\{\tau_u^1=\allowbreak t\}$ and $\{\tau_u^2=\allowbreak t\}$  are in
$\sigma(z_t)$ for all $t\ge 0$, both $\tau_u^1$ and $\tau_u^2$ belong to $\Gamma$.
Based on these stopping times, we now introduce a strategy pair induced by $u\in\mathcal{U}$, whose existence will be established in Lemma~\ref{lem:existence-strategy} below.

	 \begin{definition}\label{def:induced-strategy}
	 	For each $ u \in \mathcal{U} $, a strategy pair $((\phi_u,\allowbreak  \tau_u^1);(\psi_u,\allowbreak  \tau_u^2)) \in (\Pi^1_{\mathrm{s}} \times \Gamma) \times (\Pi^2_{\mathrm{s}} \times \Gamma)$ is called a strategy pair induced by $ u $ if
	 \begin{eqnarray*}
	 	\mathcal{T}u(i,\lambda)
	 	&=& \inf_{\psi \in \mathcal{P}(B(i))} \sum_{a\in A(i)}\sum_{b\in B(i)}\mathbb{G}u(i,\lambda;a,b)\phi_u(a|i,\lambda)\psi(b)\\
	 	&=& \sup_{\phi \in \mathcal{P}(A(i))} \sum_{a\in A(i)}\sum_{b\in B(i)}\mathbb{G}u(i,\lambda;a,b)\phi(a)\psi_u(b|i,\lambda)
	 \end{eqnarray*}
	 for	all $(i,\allowbreak \lambda)\in D^c\times  \Lambda^*$, where $ \tau_u^1 $ and $ \tau_u^2 $ are given in (\ref{tau1}) and (\ref{tau2}), respectively.
	 \end{definition}
	
	 \begin{lemma}\label{lem:existence-strategy}
	 	For any $u\in\mathcal{U}$,
	 	there is a strategy pair $((\phi_u,\allowbreak  \tau_u^1);(\psi_u,\allowbreak  \tau_u^2)) \in(\Pi^1_{\mathrm{s}}\times\Gamma)\times(\Pi^2_{\mathrm{s}}\times\Gamma)$ induced by $u$.
	 \end{lemma}
	 \begin{proof}
	 	For any fixed $u\in\mathcal{U}$, by the minimax measurable selection theorem  \cite[Lemma 4.3]{Nowak1984Onzero}, there exists an $(\phi_u,\allowbreak \psi_u)\in \Pi^1_{\mathrm{s}}\times \Pi^2_{\mathrm{s}}$ such that, for each $(i,\allowbreak \lambda)\in D^c\times  \Lambda^*$,
	 	\begin{eqnarray}
	 	\mathcal{T}u(i,\lambda)
	 		&=&\inf_{\psi \in \mathcal{P}(B(i))} \sup_{\phi \in \mathcal{P}(A(i))} \sum_{a\in A(i)}\sum_{b\in B(i)}\mathbb{G}u(i,\lambda;a,b)\phi(a)\psi(b) \nonumber\\
	 		&=& \inf_{\psi \in \mathcal{P}(B(i))} \sum_{a\in A(i)}\sum_{b\in B(i)}\mathbb{G}u(i,\lambda;a,b)\phi_u(a|i,\lambda)\psi(b)\label{phi_u}\\
	 		&=& \sup_{\phi \in \mathcal{P}(A(i))} \sum_{a\in A(i)}\sum_{b\in B(i)}\mathbb{G}u(i,\lambda;a,b)\phi(a)\psi_u(b|i,\lambda)\label{psi_u},
	 	\end{eqnarray}
	 	which, together with $\tau_u^1,\allowbreak  \tau_u^2\in\Gamma$, gives the desired result.
	 \end{proof}

	\begin{lemma}\label{lem:comparison-stopping} For any $u\in\mathcal{U}$, let $((\phi_u,\allowbreak  \tau_u^1);(\psi_u,\allowbreak  \tau_u^2)) $ be a strategy pair induced by $u$. Then, under  Assumption~\ref{ass:absorb}, the following assertions hold.
		\begin{itemize}
		\item[(a)] If $u\ge\mathcal Lu$ on $D^c\times\Lambda^*$, then, for every $(\pi^1,\allowbreak \tau^1)\in\Pi^1\times\Gamma$,
		\begin{eqnarray}\label{compare-1}
		u(i_0,\lambda_0)\ge
		\mathbb{P}_{(i_0,\lambda_0)}^{\pi^1,\psi_u}
		\bigl(X(\tau^1,\tau_u^2)\ge\lambda_0\bigr)  \ \ \  \forall (i_0,\allowbreak \lambda_0)\in D^c\times\Lambda^*.
		\end{eqnarray}
		\item[(b)] If $u\le\mathcal Lu$ on $D^c\times\Lambda^*$, then, for every $(\pi^2,\allowbreak \tau^2)\in\Pi^2\times\Gamma$,
		\begin{eqnarray}\label{compare-2}
		u(i_0,\lambda_0)\le
		\mathbb{P}_{(i_0,\lambda_0)}^{\phi_u,\pi^2}
		\bigl(X(\tau_u^1,\tau^2)\ge\lambda_0\bigr)  \ \ \  \forall (i_0,\allowbreak \lambda_0)\in D^c\times\Lambda^*.
		\end{eqnarray}
		\end{itemize}
	\end{lemma}
	\begin{proof} (a) Recall $\tau:=\allowbreak \tau_D\wedge\tau^1\wedge\tau_u^2$ in \eqref{termination time}. Then, for any $t<\tau$, by the definition of $\tau_u^2$ in \eqref{tau2},
		\begin{eqnarray}\label{LgeT}
		\mathbf{1}_{\{g^2(i_t)\ge\lambda_t\}}>\mathcal Lu(i_t,\lambda_t).
		\end{eqnarray}
		We claim that $\mathcal Lu(i_t,\allowbreak \lambda_t)\ge\mathcal Tu(i_t,\allowbreak \lambda_t)$.
		Otherwise, if
		$
		\mathcal Tu(i_t,\allowbreak \lambda_t)>\mathcal Lu(i_t,\allowbreak \lambda_t),\allowbreak
		$
		together with \eqref{LgeT}, we have
		$
		\min\left\{\mathbf{1}_{\{g^2(i_t)\ge\lambda_t\}},\allowbreak \mathcal Tu(i_t,\allowbreak \lambda_t)\right\}
		>\mathcal Lu(i_t,\allowbreak \lambda_t).
		$
		However, by the definition of $\mathcal{L}u$ in \eqref{Lu},
		\begin{eqnarray*}
		\mathcal Lu(i_t,\lambda_t)
		\ge\min\left\{\mathbf{1}_{\{g^2(i_t)\ge\lambda_t\}},\mathcal Tu(i_t,\lambda_t)\right\},
		\end{eqnarray*}
		which is a contradiction. Therefore $\mathcal Lu(i_t,\allowbreak \lambda_t)\ge\mathcal Tu(i_t,\allowbreak \lambda_t)$.
		Using the definitions of $\mathcal T$ and $\mathbb G$ in \eqref{Tu} and \eqref{G}, together with $u\ge\mathcal Lu\ge\mathcal Tu$ on $\{t<\tau\}$, by (\ref{psi_u}) we obtain
		\begin{eqnarray}\label{eq:supermartingale}
		u(i_t,\lambda_t)\mathbf1_{\{t<\tau\}}\ge\mathcal Tu(i_t,\lambda_t)\mathbf1_{\{t<\tau\}}
		\ge\mathbb E_{(i_0,\lambda_0)}^{\pi^1,\psi_u}
		\left[u(i_{t+1},\lambda_{t+1})\mathbf1_{\{t<\tau\}}| z_t\right].
		\end{eqnarray}
		For all $t\ge 0$, we define
	$
		M_t:=\allowbreak u(i_{t\wedge\tau},\allowbreak \lambda_{t\wedge\tau}).
$
		Then we get $M_{t+1}\mathbf1_{\{t<\tau\}}=\allowbreak u(i_{(t+1)\wedge\tau},\allowbreak \lambda_{(t+1)\wedge\tau})\mathbf1_{\{t<\tau\}}=\allowbreak u(i_{t+1},\allowbreak \lambda_{t+1})\mathbf1_{\{t<\tau\}}$, which, together with
		\eqref{eq:supermartingale}, gives that
		$M_t\mathbf1_{\{t<\tau\}}\ge \mathbb{E}_{(i_0,\allowbreak \lambda_0)}^{\pi^1,\allowbreak \psi_u}[M_{t+1}\mathbf1_{\{t<\tau\}}|z_t]$.
		Furthermore, since $M_t\mathbf1_{\{t\ge\tau\}}=\allowbreak u(i_{\tau},\allowbreak \lambda_{\tau})\mathbf1_{\{t\ge\tau\}}=\allowbreak M_{t+1}\mathbf1_{\{t\ge\tau\}}$, we have
		$
		M_t\ge \mathbb{E}_{(i_0,\allowbreak \lambda_0)}^{\pi^1,\allowbreak \psi_u}[M_{t+1}| z_t]$ for all $t\ge0.
		$
		Consequently,
		\begin{eqnarray}\label{eq:supermartingale2}
		\mspace{-18mu}u(i_0,\lambda_0)=M_0\ge
		 \mathbb{E}_{(i_0,\lambda_0)}^{\pi^1,\psi_u}[M_t]=\mathbb{E}_{(i_0,\lambda_0)}^{\pi^1,\psi_u}[u(i_{t\wedge\tau},\lambda_{t\wedge\tau})]\qquad \forall t\ge0.
		\end{eqnarray}
		
		Since $\tau\le\tau_D<\infty$ almost surely (by
Assumption~\ref{ass:absorb}) and $0\le u\le1$,
letting $t\to\infty$ in \eqref{eq:supermartingale2}, together with applying the dominated convergence theorem, yields
		\begin{eqnarray}\label{eq:supermartingale3}
		u(i_0,\lambda_0)\ge
		\mathbb{E}_{(i_0,\lambda_0)}^{\pi^1,\psi_u}[u(i_\tau,\lambda_\tau)].
		\end{eqnarray}
		Set
		$
		U_1:=\allowbreak \{\tau^1\le\tau_u^2,\allowbreak \ \tau^1<\tau_D\}$, $
		U_2:=\allowbreak \{\tau_u^2<\tau^1,\allowbreak \ \tau_u^2<\tau_D\}$,		$U_3:=\allowbreak \{\tau_D\le\tau^1\wedge\tau_u^2\}$, and
		\begin{eqnarray}
		\begin{aligned}
		W_1&:=\left\{\sum_{k=0}^{\tau^1-1}r(i_k,a_k,b_k)
		+g^1(i_{\tau^1})\ge\lambda_0\right\}\cap U_1,~~
		W_2:=\left\{\sum_{k=0}^{\tau_u^2-1}r(i_k,a_k,b_k)
		+g^2(i_{\tau_u^2})\ge\lambda_0\right\}\cap U_2,\\
		W_3&:=\left\{\sum_{k=0}^{\tau_D-1}r(i_k,a_k,b_k)
		\ge\lambda_0\right\}\cap U_3.
		\end{aligned}
		\end{eqnarray}
		Then,  by the definition of $\mathcal{L}$, $\tau_{u}^2$ and $u\in\mathcal{U}$, we have
		\begin{eqnarray*}
		&&u(i_{\tau^1},\lambda_{\tau^1})\mathbf{1}_{U_1}
		\ge\mathcal Lu(i_{\tau^1},\lambda_{\tau^1})\mathbf{1}_{U_1}
		\ge\mathbf{1}_{\{g^1(i_{\tau^1})\ge\lambda_{\tau^1}\}}\mathbf{1}_{U_1}
		=\mathbf{1}_{W_1},~(\text{by}~\eqref{Lu})\\
		&&u(i_{\tau_u^2},\lambda_{\tau_u^2})\mathbf{1}_{U_2}
		\ge\mathcal Lu(i_{\tau_u^2},\lambda_{\tau_u^2})\mathbf{1}_{U_2}
		\ge\mathbf{1}_{\{g^2(i_{\tau_u^2})\ge\lambda_{\tau_u^2}\}}\mathbf{1}_{U_2}
		=\mathbf{1}_{W_2},~(\text{by}~\eqref{tau2})\\
		&&u(i_{\tau_D},\lambda_{\tau_D})\mathbf{1}_{U_3}
		=\mathbf{1}_{\{\lambda_{\tau_D}\le0\}}\mathbf{1}_{U_3}
		=\mathbf{1}_{W_3}. ~~~~~~~(\text{by}~\eqref{def:U-boundary})
		\end{eqnarray*}
		Consequently,
		\begin{eqnarray*}
		\mathbb{E}_{(i_0,\lambda_0)}^{\pi^1,\psi_u}[u(i_\tau,\lambda_\tau)]
		&=&\mathbb{E}_{(i_0,\lambda_0)}^{\pi^1,\psi_u}[
		u(i_{\tau^1},\lambda_{\tau^1})\mathbf{1}_{U_1}
		+u(i_{\tau_u^2},\lambda_{\tau_u^2})\mathbf{1}_{U_2}
		+u(i_{\tau_D},\lambda_{\tau_D})\mathbf{1}_{U_3}
		]\\
		&\ge&\mathbb{P}_{(i_0,\lambda_0)}^{\pi^1,\psi_u}(W_1)
		+\mathbb{P}_{(i_0,\lambda_0)}^{\pi^1,\psi_u}(W_2)
		+\mathbb{P}_{(i_0,\lambda_0)}^{\pi^1,\psi_u}(W_3)
		=\mathbb{P}_{(i_0,\lambda_0)}^{\pi^1,\psi_u}
		\bigl(X(\tau^1,\tau_u^2)\ge\lambda_0\bigr),
		\end{eqnarray*}
		where the last equality follows from \eqref{def:total-reward}. Combining the preceding display with \eqref{eq:supermartingale3} yields
		$
		u(i_0,\lambda_0)\ge
		\mathbb{P}_{(i_0,\lambda_0)}^{\pi^1,\psi_u}
		\bigl(X(\tau^1,\tau_u^2)\ge\lambda_0\bigr).
		$
		This proves part~(a).

(b) Similarly, we can prove (b). The details are omitted.
	\end{proof}
	Lemma~\ref{lem:comparison-stopping} gives a pointwise comparison for each initial state $(i_0,\lambda_0)\in D^c\times\Lambda^*$. Taking expectation with respect to the initial distribution $\nu$ then yields the following comparison theorem.
	\begin{theorem}[Comparison Theorem]\label{thm:relationship} Suppose that Assumption \ref{ass:absorb} holds. For any $u\in\mathcal{U}$, let $((\phi_u,\allowbreak \tau^1_u),\allowbreak (\psi_u,\allowbreak \tau^2_u))$ be the   pair induced by $u$. The following statements hold.
	 \begin{itemize}
	 \item[(a)] If  $u(i,\allowbreak \lambda)\ge \mathcal{L}u(i,\allowbreak \lambda)$ for each  $(i,\allowbreak \lambda)\in D^c\times  \Lambda^*$,  then
	 \[
	 	\sum_{i\in D^c}u(i,\lambda^*)\nu(i)\ge
	 	J_{\nu}(\pi^1,\tau^1;\psi_u,\tau_u^2)~ ~\forall (\pi^1,\tau^1)\in \Pi^1\times\Gamma.
	 \]
	 \item[(b)] If  $u(i,\allowbreak \lambda)\leq \mathcal{L}u(i,\allowbreak \lambda)$ for each  $(i,\allowbreak \lambda)\in D^c\times  \Lambda^*$, then
	 \[
	 	\sum_{i\in D^c}u(i,\lambda^*)\nu(i)\leq
	 		J_{\nu}(\phi_u,\tau_u^1;\pi^2,\tau^2)~~\forall (\pi^2,\tau^2)\in \Pi^2\times\Gamma.
	 \]
 \end{itemize}
	 \end{theorem}
\begin{proof}
	(a) By Lemma~\ref{lem:comparison-stopping}, for every $(i,\allowbreak \lambda)\in D^c\times\Lambda^*$,
	 	 	\begin{eqnarray}\label{inequ-3}
	 	 		u(i,\lambda)\ge \mathbb{P}_{(i,\lambda)}^{\pi^1,\psi_u}( X(\tau^1,\tau_u^2)
	 	 		\ge \lambda).
	 	 	\end{eqnarray}
	 	 	Since $\lambda^*>0$, \eqref{def:U-boundary} gives $u(i,\allowbreak \lambda^*)=\allowbreak 0$ for every $i\in D$. Therefore, setting $\lambda=\allowbreak \lambda^*$, multiplying \eqref{inequ-3} by $\nu(i)$, and summing over $i\in D^c$ yield
	 	 	$
	 	 	\sum_{i\in S}u(i,\lambda^*)\nu(i)
	 	 	=\sum_{i\in D^c}u(i,\lambda^*)\nu(i)
	 	 	\ge
	 	 	J_{\nu}(\pi^1,\tau^1;\psi_u,\tau_u^2),
	 	 	$
	 	 	which proves  (a).\\
	 	(b) The proof is analogous, and thus the details are omitted.
	\end{proof}
	Under the general condition (i.e. Assumption \ref{ass:absorb}), we establish a comparison theorem
(Theorem~\ref{thm:relationship}), which is the key tool for proving the existence of the value and a Nash equilibrium of the game below.

 	 \section{Existence of the Value and a Nash Equilibrium}\label{sec:extended}
	In this section, we establish the existence of the value and a Nash equilibrium of the game via two monotone sequences generated by the Shapley operator $\mathcal{L}$. Their difference also provides
an online stopping rule for the algorithm developed later.
	 \begin{theorem}\label{game-Nash}
	 	Suppose Assumption \ref{ass:absorb} holds.
		Fix $\nu\in\mathcal{P}(S)$ and $\lambda^*>0$.  Let
	 	\begin{eqnarray}
	 		&&u_0^0(i,\lambda):=\mathbf{1}_{\{\lambda\le0\}},\ u_0^1(i,\lambda):=1,\ u_t^k(i,\lambda):=\mathcal{L}u_{t-1}^k(i,\lambda)\ \ \forall t\ge1,\ k=0,1, (i,\allowbreak \lambda)\in D^c\times\Lambda^*;\label{initial-two-sequence-1}\\
	 	  &&u_t^k(i,\lambda):=\mathbf{1}_{\{\lambda\le0\}}\quad \forall  t\ge0,\ k=0,1,(i,\allowbreak \lambda)\in D\times\Lambda^*. \label{initial-two-sequence-2}
	 	\end{eqnarray}
	 	Then we have the following statements.
	 	 \begin{itemize}
\item[(a)] $u_t^0\leq u_{t+1}^0$, $u_t^1\ge u_{t+1}^1$, and $u_t^0\leq u_t^1$ for all $t\ge 0$. The limits $u_*^k:=\allowbreak \lim_{t\to\infty}u_t^k$ exist, and $u_*^k$ satisfies the Shapley equation, i.e., $u_*^k=\allowbreak \mathcal{L}u_*^k$ for $k=\allowbreak 0,\allowbreak 1$. Moreover, $u_*^0=\allowbreak u_*^1$.
	 	\item[(b)] The value of the game $J^*_\nu$ exists, and $J^*_\nu=\sum_{i\in D^c}u_*(i,\lambda^*)\nu(i)$, where $u_*:=u_*^0.$
	 	\item[(c)] $((\phi_{u_*},\allowbreak  \tau_{u_*}^1),\allowbreak  (\psi_{u_*},\allowbreak  \tau_{u_*}^2)) \in (\Pi^1_{\mathrm{s}} \times \Gamma) \times (\Pi^2_{\mathrm{s}} \times \Gamma)$ induced by $u_*$ is a Nash equilibrium with respect to $(\nu,\allowbreak \lambda^*)$.
 \end{itemize}	
\end{theorem}
	 \begin{proof}
	 	(a) First, for $i\in D$ and $\lambda\in\Lambda^*$,  by \eqref{initial-two-sequence-2}, we have $u_t^0(i,\allowbreak \lambda)=\allowbreak u_t^1(i,\allowbreak \lambda)=\allowbreak \mathbf{1}_{\{\lambda\le0\}}$ for all $t\ge0$.
	 	For $i\in D^c$ and $\lambda\le0$, by \eqref{initial-two-sequence-1} and the definition of $\mathcal{L}$, we get $u_t^0(i,\allowbreak \lambda)=\allowbreak u_t^1(i,\allowbreak \lambda)=\allowbreak 1$ for all $t\ge0$.
	 	Thus the part (a) holds trivially for the above two cases.
	 	It remains to consider $i\in D^c$ and $\lambda>0$,
for which we consider the sequences $\{u_t^0\}_{t\ge0}$
and $\{u_t^1\}_{t\ge0}$ separately.\\
	 	For the sequence $\{u_t^0\}_{t\ge0}$, we have $u_1^0(i,\allowbreak \lambda)=\allowbreak \mathcal{L}u_0^0(i,\allowbreak \lambda)\ge0=\allowbreak u_0^0(i,\allowbreak \lambda)$. The monotonicity of $\mathcal{L}$ implies that $\{u_t^0\}$ is a nondecreasing sequence and bounded above by $1$, and hence converges pointwise to some $u_*^0\in\mathcal{U}$. Then we have $\mathcal{L}u_*^0\ge \mathcal{L}u_{t}^0=\allowbreak u_{t+1}^0$ for all $t\ge 0$, which implies $\mathcal{L}u_*^0\ge u_*^0$.\\
		Next we prove that $u_*^0\ge \mathcal{L}u_*^0$.
	 	First, for $i\in D^c$, $\lambda>0$, let $
	 		G_1(i,\allowbreak \lambda):=\allowbreak \mathbf{1}_{\{g^1(i)\ge\lambda\}}$ and  $G_2(i,\allowbreak \lambda):=\allowbreak \mathbf{1}_{\{g^2(i)\ge\lambda\}}.
	    $
	 	Then by the definition of
	 	 $\mathcal{L}$ in \eqref{Lu},
	 	 it remains to prove
	 	 \begin{eqnarray}\label{remain-u_*^0}
	 	 \mspace{-20mu}u_*^0(i,\lambda)\ge G_1(i,\lambda),  u_*^0(i,\lambda)\ge\min\{G_2(i,\lambda),\mathcal{T}u_*^0(i,\lambda)\}.
	 	 \end{eqnarray}
	 	(i) If $G_1(i,\allowbreak \lambda)=\allowbreak 0$, it is trivial that $u_*^0(i,\allowbreak \lambda)\ge G_1(i,\allowbreak \lambda)$. \\

	 	(ii) If $G_1(i,\allowbreak \lambda)=\allowbreak 1$, then $g^1(i)\ge\lambda$, which implies
	 	$u_1^0(i,\allowbreak \lambda)=\allowbreak \mathcal{L}u_0^0(i,\allowbreak \lambda)=\allowbreak \max\{1,\allowbreak \min\{G_2(i,\allowbreak \lambda),\allowbreak \mathcal{T}u_0^0(i,\allowbreak \lambda)\}\}=\allowbreak 1$. Since $\{u_t^0\}$ is a nondecreasing sequence, for all
	 	$t\ge1$, $u_t^0(i,\allowbreak \lambda)\ge u_1^0(i,\allowbreak \lambda)=\allowbreak 1$, hence $u_*^0(i,\allowbreak \lambda)=\allowbreak 1\ge G_1(i,\allowbreak \lambda)$.\\

	 	(iii) If $G_2(i,\allowbreak \lambda)=\allowbreak 0$, $u_*^0(i,\allowbreak \lambda)\ge \min\{G_2(i,\allowbreak \lambda),\allowbreak \mathcal{T}u_*^0(i,\allowbreak \lambda)\}$.\\

	 	(iv) If $G_2(i,\allowbreak \lambda)=\allowbreak 1$, then, for any fixed $\mu\in\mathcal{P}(A(i))$, by the definition of $\mathcal{L}$ and $\mathcal{T}$, there exists $\psi_t^\mu\in\mathcal{P}(B(i))$ such that
	 	\begin{eqnarray}\label{thm1-1}
	 		u_{t+1}^0(i,\lambda)=\mathcal{L}u_t^0(i,\lambda)\ge \mathcal{T}u_t^0(i,\lambda)
	 		\ge\sum_{a\in A(i)}\sum_{b\in B(i)}\mathbb{G}u_t^0(i,\lambda;a,b)\mu(a)\psi_t^\mu(b)~~\forall t\ge0.
	 	\end{eqnarray}
	 	By the compactness of $\mathcal{P}(B(i))$, passing to a subsequence if necessary, we may assume that $\psi_t^\mu\to\psi^{\mu,\allowbreak *}\in\mathcal{P}(B(i))$.
		 Since $A(i)$ and $B(i)$ are finite, by \eqref{thm1-1} and the dominated convergence theorem we have
		\begin{eqnarray*}
	 		u_*^0(i,\lambda)\ge\sum_{a\in A(i)}\sum_{b\in B(i)}\mathbb{G}u_*^0(i,\lambda;a,b)\mu(a)\psi^{\mu,*}(b)
	 		\ge\inf_{\psi\in\cP(B(i))}\sum_{a\in A(i)}\sum_{b\in B(i)}\mathbb{G}u_*^0(i,\lambda;a,b)\mu(a)\psi(b),
	 	\end{eqnarray*}
	 	which, together with the arbitrariness of $\mu\in\cP(A(i))$, implies that $u_*^0(i,\allowbreak \lambda)\ge\mathcal{T}u_*^0(i,\allowbreak \lambda)$.
	 	Since $G_2(i,\allowbreak \lambda)=\allowbreak 1$, we have $u_*^0(i,\allowbreak \lambda)\ge\min\{G_2(i,\allowbreak \lambda),\allowbreak \mathcal{T}u_*^0(i,\allowbreak \lambda)\}$.
	 	Based on the above cases (i)-(iv), we can get \eqref{remain-u_*^0} for $i\in D^c$ and $\lambda>0$, and hence $u_*^0\ge\mathcal{L}u_*^0$. Since the reverse inequality was proved above, $u_*^0=\allowbreak \mathcal{L}u_*^0$.\\
	 	Similarly, for $i\in D^c$ and $\lambda>0$, \eqref{initial-two-sequence-1} gives $u_0^1(i,\allowbreak \lambda)=\allowbreak 1$, and hence $u_1^1(i,\allowbreak \lambda)=\allowbreak \mathcal{L}u_0^1(i,\allowbreak \lambda)\le 1=\allowbreak u_0^1(i,\allowbreak \lambda)$. The monotonicity of $\mathcal{L}$ then implies that  $\{u_t^1,\allowbreak t\ge0\}$ is a nonincreasing sequence and bounded below by $0$, and hence converges  pointwise to some $u_*^1\in\mathcal{U}$.
		Applying the same case analysis to the nonincreasing sequence $\{u_t^1\}$ (with the inequalities reversed) yields $u_*^1=\allowbreak \mathcal{L}u_*^1$.
	 	By induction,  $u_t^0\le u_t^1$ for every $t$, and so $u_*^0\le u_*^1$. To prove the reverse inequality (i.e.$u_*^0\ge u_*^1)$, we choose the strategy pairs $(\phi_{u_*^1},\allowbreak \tau_{u_*^1}^1;\cdot,\allowbreak \cdot)$ and $(\cdot,\allowbreak \cdot;\psi_{u_*^0},\allowbreak \tau_{u_*^0}^2)$ induced by $u_*^1$ and $u_*^0$, respectively.
		By Lemma \ref{lem:comparison-stopping}(b) and $u_*^1=\allowbreak \mathcal{L}u_*^1$, taking  $u=\allowbreak u_*^1$ and
		$
		(\pi^2,\allowbreak \tau^2)=\allowbreak (\psi_{u_*^0},\allowbreak \tau_{u_*^0}^2)
		$ in \eqref{compare-2}, we obtain
	 	\begin{eqnarray}\label{use-compare-2}
	 	\mspace{-19mu}u_*^1(i,\lambda)
	 	\le \mathbb{P}_{(i,\lambda)}^{\phi_{u_*^1},\psi_{u_*^0}}
	 	\bigl(X(\tau_{u_*^1}^1,\tau_{u_*^0}^2)\ge\lambda\bigr)~\forall~ i\in D^c, \lambda>0.
	 	\end{eqnarray}
		By Lemma~\ref{lem:comparison-stopping}(a) and $u_*^0=\allowbreak \mathcal{L}u_*^0$,  taking $u=\allowbreak u_*^0$ and
		$
		(\pi^1,\allowbreak \tau^1)=\allowbreak (\phi_{u_*^1},\allowbreak \tau_{u_*^1}^1)
		$ in \eqref{compare-1},  we obtain
		\begin{eqnarray}\label{use-compare-1}
		\mspace{-19mu}\mathbb{P}_{(i,\lambda)}^{\phi_{u_*^1},\psi_{u_*^0}}
		\bigl(X(\tau_{u_*^1}^1,\tau_{u_*^0}^2)\ge\lambda\bigr)
		\le u_*^0(i,\lambda)~\forall~ i\in D^c, \lambda>0.
		\end{eqnarray}
		which, together with  \eqref{use-compare-2}, yields
		$u_*^1(i,\allowbreak \lambda)\le u_*^0(i,\allowbreak \lambda)$ for all~ $i\in D^c$ and $\lambda>0$.  Hence $u_*^0=\allowbreak u_*^1$.
		
	 	We prove (b) and (c) together. Since $u_*=\allowbreak \mathcal{L}u_*$,
	 	Lemma~\ref{lem:existence-strategy} yields a strategy pair
	 	$(\phi_{u_*},\allowbreak \tau_{u_*}^1;\psi_{u_*},\allowbreak \tau_{u_*}^2)$
	 	induced by $u_*$. For any $(\pi^1,\allowbreak \tau^1)\in\Pi^1\times\Gamma$
	 	and $(\pi^2,\allowbreak \tau^2)\in\Pi^2\times\Gamma$, applying
	 	Theorem~\ref{thm:relationship} (a) and (b), respectively, gives
	 	\begin{eqnarray}\label{Nash}
	 	J_\nu(\pi^1,\tau^1;\psi_{u_*},\tau_{u_*}^2)
	 		\leq \sum_{i\in D^c}u_*(i,\lambda^*)\nu(i)
	 		\leq J_\nu(\phi_{u_*},\tau_{u_*}^1;\pi^2,\tau^2).
	 	\end{eqnarray}
	 	It follows that
	 	$
	 	\overline{J}_\nu\leq
	 	\sum_{i\in D^c}u_*(i,\allowbreak \lambda^*)\nu(i)
	 	\leq \underline{J}_\nu.
	 	$
	 	Since $\underline{J}_\nu\leq\overline{J}_\nu$, the value exists and
	 	$J^*_\nu=\allowbreak \sum_{i\in D^c}u_*(i,\allowbreak \lambda^*)\nu(i)$. This proves (b),
	 	and (c) follows directly from \eqref{Nash}.
	 \end{proof}
	\begin{remark}
(a) Theorem~\ref{game-Nash} shows that the two monotone sequences
converge to a unique solution to the Shapley equation, which is used to construct the value
and induce a Nash equilibrium of the game.
The two sequences provide lower and upper approximations to the solution to the Shapley equation,
and the  difference of the two sequences will be used to provide an online stopping rule for the
algorithm developed below.

(b) If neither player stops, i.e.,
$\tau^1=\allowbreak \tau^2=\allowbreak \infty$, then the payoff reduces to
$
\mathbb{P}_{(i,\allowbreak \lambda)}^{\pi^1,\allowbreak \pi^2}
\left(
\sum_{t=\allowbreak 0}^{\tau_D-1}r(i_t,\allowbreak a_t,\allowbreak b_t)\ge\lambda
\right),
$
and the corresponding Shapley operator is $\mathcal T$.
Thus, the above analysis recovers the existence of the value and
stationary saddle-point policies for the   game without
stopping \cite[Theorem~3.6]{HGP-2017}.
\end{remark}
\section{Truncated Games and Finite-Approximation}

In the previous section, we established a first main result:  The two monotone sequences $\{u_t^0\}$ and $\{u_t^1\}$ in Theorem~\ref{game-Nash}  converge  to the unique  solution $u_*$ of the Shapley equation, which is used to construct the value of the game  by $J_\nu^*=\allowbreak \sum_{i\in D^c}u_*(i,\allowbreak \lambda^*)\nu(i)$.

However, the iterative construction of  $\{u_t^0(i,\lambda)\}$ or $\{u_t^1(i,\lambda)\}$   requires evaluating them at every $i\in D^c$ and  $\lambda\in \Lambda^*$. When the state space $S$ is countably infinite, direct implementation is impossible because one cannot enumerate infinitely many states. 
To overcome this difficulty, we introduce a truncation method,by which we show that  the original infinite-state game can be approximated by a sequence of truncated games based on sublevel sets of the Lyapunov function. Moreover,  we first obtain explicit error bounds between the values of the original game  and any truncated game. These bounds yield the convergence of the truncated game values
to the original game value; see Theorem~\ref{thm:truncation} below.

Throughout this section we work under the Lyapunov condition \eqref{eq:drift_section} of Proposition~\ref{prop:attenuation}.

\subsection[The Model of the Truncated Game]{The Model of the Truncated Game ${\mathcal{G}}_n$}
For each   $n\ge1$, define the sublevel set of the Lyapunov function $w$ on $D^c$ at truncation level $n$ by:
\begin{equation}\label{sublevel-set}
D_n^c \;:=\; \{\,i\in D^c : w(i)\le n\,\}.
\end{equation}
and set $S_n:=\allowbreak D_n^c\cup D$.
We now construct a truncated game $\mathcal{G}_n$ by the tuple
\[
{\mathcal{G}}_n := \{S_n^*,\; D^*,\; (A_n(i), B_n(i))_{i\in S_n^*},\; p_n,\; r_n,\; g^1_n,\; g^2_n,\; \lambda^*\},
\]
where:
\begin{itemize}
	\item Let  $i^*$ be an auxiliary state, and define
	$
	S_n^*:=S_n\cup\{i^*\},$ $
	D^*:=D\cup\{i^*\}.
	$

	\item Let $a_*^1$ and $a_*^2$ be auxiliary actions, and set
	$
	A^*:=A\cup\{a_*^1\}$, $B^*:=B\cup\{a_*^2\}.
	$
	For $i\in S_n^*$, define the admissible action sets $A_n(i)\subset A^*$ and $B_n(i)\subset B^*$ by
	\[
	A_n(i):=
	\begin{cases}
	A(i), & i\in S_n,\\
	\{a_*^1\}, & i=i^*,
	\end{cases}
	~~~
	B_n(i):=
	\begin{cases}
	B(i), & i\in S_n,\\
	\{a_*^2\}, & i=i^*.
	\end{cases}
	\]
	\item Let ${p}_n(i^*| i^*,\allowbreak a^1_*,\allowbreak a^2_*):=\allowbreak 1$, and for any $i\in S_n$, $a\in {A}_n(i),\allowbreak  b\in {B}_n(i)$,  let
	\begin{eqnarray}\label{def:p_n}
		p_n(j|i,a,b) :=
		\begin{cases}
			p(j|i,a,b), &j\in S_n,\\
			\displaystyle\sum_{s\notin S_n} p(s|i,a,b),  &j = i^*.
		\end{cases}
	\end{eqnarray}
	\item Set $r_n(i^*,\allowbreak a^1_*,\allowbreak a^2_*)=\allowbreak g^1_n(i^*)=\allowbreak g^2_n(i^*)=\allowbreak 0$, and for $i\in S_n$,
	\begin{eqnarray}\label{eq:truncated-rewards}
	r_n(i,\cdot,\cdot):=r(i,\cdot,\cdot),
	g^1_n(i):=g^1(i),
	g^2_n(i):=g^2(i).
	\end{eqnarray}
\end{itemize}
Let
$
\Lambda_n^*:=\allowbreak
\{\lambda^*-\sum_{s=\allowbreak 0}^{m-1}r_n(i_s,\allowbreak a_s,\allowbreak b_s):
m\ge0,\allowbreak \ (i_s,\allowbreak a_s,\allowbreak b_s)\in\mathbb{K}_n,\allowbreak  i_s\in D_n^c,\allowbreak  s=\allowbreak 0,\allowbreak 1,\allowbreak \ldots,\allowbreak m-1\},\allowbreak
$
where $\lambda^*>0$ is the initial threshold and
$
\mathbb{K}_n:=\allowbreak
\left\{(i,\allowbreak a,\allowbreak b):i\in S_n^*,\allowbreak \ a\in A_n(i),\allowbreak \ b\in B_n(i)\right\}.
$

Fix $n$. The corresponding objects for the truncated game $\mathcal G_n$ are denoted with a hat, such as the canonical path space
$(\widehat\Omega,\allowbreak \widehat{\mathcal F})$, the history sets
$\widehat H_t$, and the coordinate and history processes
$(\widehat i_t,\allowbreak \widehat\lambda_t,\allowbreak \widehat a_t,\allowbreak \widehat b_t,\allowbreak \widehat z_t)$.
 Moreover, to indicate the dependence on the truncation level $n$, the subscript $n$ is added to some of them, such as the stopping-time and policy classes $\Gamma_n$, $\Pi_n^1$, $\Pi_n^2$, $\Pi_{n,\allowbreak \mathrm{s}}^1$, and $\Pi_{n,\allowbreak \mathrm{s}}^2$. Then we define the first-passage time to $D^*$ by
$
\widehat{\tau}_{D^*}:=\allowbreak \inf\{t\ge0:\widehat i_t\in D^*\},\allowbreak
$
and for $\widehat{\tau}^1,\allowbreak  \widehat{\tau}^2\in\Gamma_n$, write the termination time by
$
\widehat{\tau}_*:=\allowbreak \widehat{\tau}_{D^*}\wedge\widehat{\tau}^1\wedge\widehat{\tau}^2,\allowbreak
$
and the total reward $\widehat X$ accumulated up to termination by
\begin{eqnarray}\label{def:truncated-total-reward}
\widehat{X}(\widehat{\tau}^1,\widehat{\tau}^2)
:=
\sum_{t=0}^{\widehat{\tau}_*-1} r_n(\widehat{i}_t,\widehat{a}_t,\widehat{b}_t)
+
\mathbf{1}_{\{\widehat{\tau}_* < \widehat{\tau}_{D^*}\}}\times
\left(
\mathbf{1}_{\{\widehat{\tau}^1 \le \widehat{\tau}^2\}}g_n^1(\widehat{i}_{\widehat{\tau}_*})
+
\mathbf{1}_{\{\widehat{\tau}^2 < \widehat{\tau}^1\}}g_n^2(\widehat{i}_{\widehat{\tau}_*})
\right).
\end{eqnarray}

Fix an initial distribution $\nu\in\mathcal{P}(S)$ of the original game $\mathcal{G}$, and define the
corresponding distribution $\nu_n\in\mathcal{P}( S_n^*)$  of the truncated game $\mathcal{G}_n$ by
$
\nu_n(i):=\allowbreak \nu(i)$ for $i\in S_n$ and $
\nu_n(i^*):=\allowbreak \nu(S\setminus S_n).
$
For any strategy pair
$((\widehat{\pi}^1,\allowbreak \widehat{\tau}^1),\allowbreak (\widehat{\pi}^2,\allowbreak \widehat{\tau}^2))$ of $\mathcal{G}_n$, the Ionescu--Tulcea
theorem yields a unique probability space $(\widehat\Omega,\allowbreak \widehat{\mathcal F},\allowbreak  \mathbb P_{n,\allowbreak \nu_n}^{\widehat{\pi}^1,\allowbreak \widehat{\pi}^2})$. Then, we  define
\[
J_{n,\nu_n}(\widehat{\pi}^1,\widehat{\tau}^1;\widehat{\pi}^2,\widehat{\tau}^2)
:=
\mathbb{P}_{n,\nu_n\otimes\delta_{\lambda^*}}^{\widehat{\pi}^1,\widehat{\pi}^2}
\big(\widehat{X}(\widehat{\tau}^1,\widehat{\tau}^2)\ge\lambda^*\big).
\]
Define $\underline J_{n,\allowbreak \nu_n}$, $\overline J_{n,\allowbreak \nu_n}$, and
$J_{n,\allowbreak \nu_n}^*$, which are analogous to those for the original game.Then, we define the
corresponding $\epsilon$-approximation value and $\epsilon$-Nash equilibria of the truncated game $\mathcal{G}_n$.

Moreover, define $\mathcal{U}_n$, $\mathbb{G}_n$, $\mathcal{T}_n$, $\mathcal{L}_n$ on $ S_n^* \times \Lambda_n^*$ by replacing the $S,\allowbreak  D^c,\allowbreak  A(i),\allowbreak  B(i),\allowbreak  p,\allowbreak  r,\allowbreak  g^1,\allowbreak  g^2,\allowbreak $ and $\Lambda^*$ by $S_n^*,\allowbreak  D_n^c,\allowbreak  A_n(i),\allowbreak  B_n(i),\allowbreak  p_n,\allowbreak  r_n,\allowbreak  g_n^1,\allowbreak  g_n^2,\allowbreak $ and $\Lambda_n^*$, respectively.

A natural question is whether the value and a Nash equilibrium actually exist for the truncated game $\mathcal{G}_n$.
 By Theorem~\ref{game-Nash}, this question is answered positively once Assumption~\ref{ass:absorb} holds for $\mathcal{G}_n$.
 In view of Proposition~\ref{prop:attenuation}, it suffices to verify the Lyapunov condition for $\mathcal{G}_n$.
 Since the transition kernel is modified under truncation, this condition has to be checked for $\mathcal{G}_n$.
 The verification and its consequences are summarized below.

\begin{corollary}\label{cor:truncated-game-Nash}
	Suppose the drift condition \eqref{eq:drift_section} holds.
	 For any fixed $n\in\mathbb{N}$,
 let $u_{n,\allowbreak 0}^0(i,\allowbreak \lambda):=\allowbreak \mathbf{1}_{\{\lambda\le0\}}$
 and $u_{n,\allowbreak 0}^1(i,\allowbreak \lambda):=\allowbreak 1$ on $D_n^c\times\Lambda_n^*$, and define $u_{n,\allowbreak t}^k:=\allowbreak \mathcal{L}_n u_{n,\allowbreak t-1}^k$ for $t\ge1$ and $k=\allowbreak 0,\allowbreak 1$.
Set
	 	$u_{n,\allowbreak t}^k(i,\allowbreak \lambda):=\allowbreak \mathbf{1}_{\{\lambda\le0\}}$ for all
	    $(i,\allowbreak \lambda)\in D^*\times\Lambda_n^*$, $t\ge0$ and $k=\allowbreak 0,\allowbreak 1$.
Then the following statements hold.
	\begin{itemize}
		\item[(a)] $u_{n,\allowbreak t}^0\leq u_{n,\allowbreak t}^1$ for all $t\geq 1$, and the limits $u^k_{n,\allowbreak *}:=\allowbreak \lim_{t\to\infty}u_{n,\allowbreak t}^k$ exist, and $u^k_{n,\allowbreak *}$ satisfies the Shapley equation,
		i.e., $u^k_{n,\allowbreak *}=\allowbreak \mathcal{L}_nu^k_{n,\allowbreak *}$ for $k=\allowbreak 0,\allowbreak 1$. Moreover, $u^0_{n,\allowbreak *}=\allowbreak u^1_{n,\allowbreak *}$.
		\item[(b)] Let $u_n^*:=\allowbreak u^0_{n,\allowbreak *}=\allowbreak u^1_{n,\allowbreak *}$. The value of the truncated game $\mathcal{G}_n$ exists and satisfies
		\[
		J_{n,\nu_n}^*=\sum_{i\in D_n^c}u_n^*(i,\lambda^*)\nu_n(i).
		\]
		\item[(c)] The $(\phi_{u_n^*},\allowbreak  \tau_{u_n^*}^1;\psi_{u_n^*},\allowbreak  \tau_{u_n^*}^2) \in (\Pi^1_{n,\allowbreak \mathrm{s}} \times \Gamma_n) \times (\Pi^2_{n,\allowbreak \mathrm{s}} \times \Gamma_n)$ induced by $u_n^*$ is a Nash equilibrium of $\mathcal{G}_n$.
	\end{itemize}
\end{corollary}
\begin{proof}
	By the definition \eqref{def:p_n} of the truncated kernel $p_n$, for every $i\in D_n^c,\allowbreak  a\in A_n(i),\allowbreak $ and $b\in B_n(i)$,
	\begin{eqnarray*}
	\sum_{j\in D_n^c} p_n(j| i,a,b)\,w(j)
	= \sum_{j\in D_n^c} p(j| i,a,b)\,w(j)
	\le \sum_{j\in D^c} p(j| i,a,b)\,w(j)
	\le \rho\,w(i).
	\end{eqnarray*}
	Thus the Lyapunov condition \eqref{eq:drift_section} holds for $\mathcal{G}_n$. By Proposition~\ref{prop:attenuation} and  Theorem~\ref{game-Nash}, we see the conclusions in this corollary are all true.
\end{proof}

\subsection[Difference between the Values of the Original and Truncated Games]{Difference between the Values of the Original Game $\mathcal{G}$ and the Truncated Game $\mathcal{G}_n$}

Having established the value and a Nash equilibrium for each truncated game $\mathcal G_n$, we next quantify how well $\mathcal G_n$ approximates the original game $\mathcal G$ by constructing two maps that transfer strategy pairs between the two games; the resulting error bounds are given in Theorem~\ref{thm:truncation} below.

For each $n\ge1$, define the first-exit times from $S_n$
for $\mathcal G$ and $\mathcal G_n$, respectively, by
\begin{eqnarray*}
\sigma_n(\omega)
&:=& \inf\{t\ge0:i_t(\omega)\notin S_n\},
\qquad \omega\in\Omega,\\
\widehat\sigma_n(\hat\omega)
&:=& \inf\{t\ge0:\widehat i_t(\hat\omega)\notin S_n\},
\qquad \hat\omega\in\widehat\Omega.
\end{eqnarray*}
In $\mathcal G_n$, $\widehat\sigma_n$ is also the
first-passage time of $i^*$. For $\widehat\tau^1,\widehat\tau^2\in\Gamma_n$, let
\[
\widehat\tau_D:=\inf\{t\ge0:\widehat i_t\in D\},
\qquad
\widehat\tau:=
\widehat\tau_D\wedge\widehat\tau^1\wedge\widehat\tau^2.
\]
Then the termination time of $\mathcal G_n$ satisfies
$\widehat\tau_*=\widehat\tau\wedge\widehat\sigma_n$.
We compare the payoffs on
$\{\widehat\tau<\widehat\sigma_n\}$
and bound the probability of
$\{\widehat\sigma_n\le\widehat\tau\}$.

For each $n\ge1$, we introduce two maps between the
strategy spaces of the original and truncated games
to compare payoffs under arbitrary strategies:
\[
\begin{aligned}
&\Phi:\;(\Pi^1\times\Gamma)\times(\Pi^2\times\Gamma)
\longrightarrow
(\Pi_n^1\times\Gamma_n)\times(\Pi_n^2\times\Gamma_n),\\
&\Psi:\;(\Pi_n^1\times\Gamma_n)\times(\Pi_n^2\times\Gamma_n)
\longrightarrow
(\Pi^1\times\Gamma)\times(\Pi^2\times\Gamma).
\end{aligned}
\]

\paragraph{(i) The definition of $\Phi$} For each $(\pi^1,\allowbreak \tau^1;\pi^2,\allowbreak \tau^2)\in(\Pi^1\times\Gamma)\times(\Pi^2\times\Gamma)$, let  $\Phi(\pi^1,\allowbreak \tau^1;\pi^2,\allowbreak \tau^2):=\allowbreak (\Phi_1(\pi^1),\allowbreak \Phi_1(\tau^1);\Phi_2(\pi^2),\allowbreak \Phi_2(\tau^2)),\allowbreak $ which are given below.\\
(i-1) For each $t\ge0$ and each truncated history
$
\widehat h_t=\allowbreak (\widehat i_0,\allowbreak \widehat\lambda_0,\allowbreak \widehat a_0,\allowbreak \widehat b_0,\allowbreak \ldots,\allowbreak \widehat i_t,\allowbreak \widehat\lambda_t)\in\widehat H_t,\allowbreak
$
if $\widehat i_s\in S_n$ for all $s\le t$,
then $\widehat h_t$ is also a feasible history of the original game,
 that is, $\widehat h_t\in H_t$.
Thus, for $k=\allowbreak 1,\allowbreak 2$, define
\begin{eqnarray}\label{eq:Phi-control}
\mspace{-18mu}\Phi_k(\pi^k)_t(\cdot|\widehat h_t)
=
	\begin{cases}
\pi^k_t(\cdot| \widehat h_t), & \text{if }\widehat i_s\in S_n\text{ for all }s\le t,\\
\delta_{a_*^k}(\cdot), & \text{otherwise}.
	\end{cases}
\end{eqnarray}
and $(\Phi_k(\pi^k)_t,\allowbreak  t=\allowbreak 0,\allowbreak 1,\allowbreak \ldots)=\allowbreak :\Phi_k(\pi^k)$.\\
(i-2) Fix $\hat\omega\in\widehat\Omega$ and $k=\allowbreak 1,\allowbreak 2$.
If $\widehat{\sigma}_n(\hat\omega)=\allowbreak \infty$,
 then $\hat\omega$ never leaves $S_n$ and
thus $\hat\omega$ is also a feasible path of the original game, that is, $\hat\omega\in\Omega$.
 Set
$\Phi_k(\tau^k)(\hat\omega):=\allowbreak \tau^k(\hat\omega)$.
If $\widehat{\sigma}_n(\hat\omega)<\infty$, let $m^*:=\allowbreak \widehat{\sigma}_n(\hat\omega)-1$ and let
$
\Omega_{\hat\omega}
:=\{\omega\in\Omega:\  i_s(\omega)=\widehat i_s(\hat\omega),\
\lambda_s(\omega)=\widehat\lambda_s(\hat\omega), a_s(\omega)=\widehat a_s(\hat\omega),\
b_s(\omega)=\widehat b_s(\hat\omega),\quad \forall s\leq m^*\}.$
If there exists $\omega\in\Omega_{\hat\omega}$ such that
$\tau^k(\omega)<\widehat{\sigma}_n(\hat\omega)$,
then for any two such paths $\omega,\omega'\in \Omega_{\hat\omega}$, by the stopping-time property of $\tau^k$, we see that $\tau^k(\omega)=\tau^k(\omega')$. Choose one such path $\omega\in \Omega_{\hat\omega}$ and denote it by $\omega_{\hat\omega}^k$, and then it is legal to define
\begin{eqnarray}\label{eq:Phi-stopping}
\Phi_k(\tau^k)(\hat\omega)
=
\begin{cases}
\tau^k(\hat\omega), & \text{if } \widehat{\sigma}_n(\hat\omega)=\infty,\\[4pt]
\tau^k(\omega_{\hat\omega}^k), & \text{if }
\widehat{\sigma}_n(\hat\omega)<\infty\text{ and such }\omega_{\hat\omega}^k\text{ exists},\\[4pt]
\infty, & \text{otherwise}.
\end{cases}
\end{eqnarray}
By construction $\{\Phi_k(\tau^k)=\allowbreak t\}\in\sigma(\widehat z_t)$, so $\Phi_k(\tau^k)\in\Gamma_n$.

\paragraph{(ii) The definition of $\Psi$} For  each
$(\widehat{\pi}^1,\allowbreak \widehat{\tau}^1;\widehat{\pi}^2,\allowbreak \widehat{\tau}^2)
\in(\Pi^1_n\times\Gamma_n)\times(\Pi^2_n\times\Gamma_n)$ for the truncated game $\mathcal G_n$,
we define $\Psi(\widehat{\pi}^1,\allowbreak \widehat{\tau}^1;\widehat{\pi}^2,\allowbreak \widehat{\tau}^2):=\allowbreak (\Psi_1(\widehat{\pi}^1),\allowbreak \Psi_1(\widehat{\tau}^1);\Psi_2(\widehat{\pi}^2),\allowbreak \Psi_2(\widehat{\tau}^2))$  as follows.
\\(ii-1)
For each $i\in S$, fix admissible  actions
$a^1(i)\in A(i)$ and $a^2(i)\in B(i)$.
For each $t\ge0$ and each original history
$h_t=\allowbreak (i_0,\allowbreak \lambda_0,\allowbreak a_0,\allowbreak b_0,\allowbreak \ldots,\allowbreak i_t,\allowbreak \lambda_t)\in H_t$,
if $i_s\in S_n$ for all $s\le t$, then $h_t$ is also a feasible history of the truncated game, that is, $h_t\in\widehat H_t$.
For $k=\allowbreak 1,\allowbreak 2$,
we define
\begin{eqnarray}\label{eq:Psi-control}
\mspace{-18mu}\Psi_k(\widehat{\pi}^k)_t(\cdot| h_t)
=
\begin{cases}
\widehat{\pi}^k_t(\cdot| h_t), & \text{if }i_s\in S_n\text{ for all }s\le t,\\
\delta_{a^k(i_t)}(\cdot), & \text{otherwise}.
\end{cases}
\end{eqnarray}
(ii-2) Fix $\omega\in\Omega$ and $k=\allowbreak 1,\allowbreak 2$.
If $\sigma_n(\omega)=\allowbreak \infty$, then $\omega$ never leaves $S_n$ and
thus $\omega$ is also a feasible path of the truncated game, that is, $\omega\in\widehat\Omega$.
Set $\Psi_k(\widehat{\tau}^k)(\omega):=\allowbreak \widehat{\tau}^k(\omega)$.
If $\sigma_n(\omega)<\infty$, let $l^*:=\allowbreak \sigma_n(\omega)-1$ and let
$
\widehat\Omega_{\omega}
:=\{\hat\omega\in\widehat\Omega:\  \widehat i_s(\hat\omega)=i_s(\omega),\
\widehat\lambda_s(\hat\omega)=\lambda_s(\omega),
 \widehat a_s(\hat\omega)=a_s(\omega),\
\widehat b_s(\hat\omega)=b_s(\omega),\quad \forall s\le l^*\}.$
If there exists $\hat\omega\in\widehat\Omega_{\omega}$ such that
$\widehat{\tau}^k(\hat\omega)<\sigma_n(\omega)$,
then any two such paths have the same history up to time $l^*$,
hence also up to $\widehat{\tau}^k$; by the stopping-time property of $\widehat{\tau}^k$, the value $\widehat{\tau}^k(\hat\omega)$ is the same for all of them. Choose one such path and denote it by $\hat\omega_{\omega}^k$. Define
\begin{eqnarray}\label{eq:Psi-stopping}
\Psi_k(\widehat{\tau}^k)(\omega)
=
\begin{cases}
\widehat{\tau}^k(\omega), & \text{if } \sigma_n(\omega)=\infty,\\[4pt]
\widehat{\tau}^k(\hat\omega_{\omega}^k), & \text{if }
\sigma_n(\omega)<\infty\text{ and such }\hat\omega_{\omega}^k\text{ exists},\\[4pt]
\infty, & \text{otherwise}.
\end{cases}
\end{eqnarray}
Note that $\{\Psi_k(\widehat{\tau}^k)=\allowbreak t\}\in\sigma(z_t)$, so $\Psi_k(\widehat{\tau}^k)\in\Gamma$.

Furthermore, by construction, the above maps satisfy
the following composition identities:
\begin{eqnarray}\label{eq:PhiPsi}
\Phi_k(\Psi_k(\Phi_k))=\Phi_k,\quad
\Psi_k(\Phi_k(\Psi_k))=\Psi_k,\quad k=1,2.
\end{eqnarray}

The following proposition provides two estimates on the events
$\{\widehat\tau<\widehat\sigma_n\}$ and $\{\widehat\sigma_n\le\widehat\tau\}$, respectively.
We recall Doob's maximal inequality for nonnegative
supermartingales, also known as Ville's inequality.
See \citet{Doob1953} for background on martingale
inequalities, and \citet[p.~261]{Durrett2019}
for the result used below.

\begin{lemma}\label{lem:ville}
Let $\{M_t\}_{t\in\mathbb N}$ be a nonnegative
supermartingale on a filtered probability space
with probability measure $\mathbb P$.
Then, for every $C>0$,
\[
\mathbb P\left(\sup_{t\in\mathbb N}M_t\ge C\right)
\le \frac{\mathbb E_{\mathbb P}[M_0]}{C}.
\]
\end{lemma}

\begin{proposition}\label{prop:two-estimates}
Fix $n\ge1$ and $(i,\allowbreak \lambda)\in D_n^c\times\Lambda_n^*$.
\begin{itemize}
\item[(a)] Let $(\pi^1,\allowbreak \tau^1;\pi^2,\allowbreak \tau^2)$ be a strategy pair of $\mathcal{G}$, and
$(\widehat{\pi}^1,\allowbreak \widehat{\tau}^1;\widehat{\pi}^2,\allowbreak \widehat{\tau}^2)
=\allowbreak \Phi(\pi^1,\allowbreak \tau^1;\pi^2,\allowbreak \tau^2)$. Then
\begin{eqnarray}\label{eq:agg-noexit}
\mathbb{P}_{n,(i,\lambda)}^{\widehat{\pi}^1,\widehat{\pi}^2}
\left(\widehat{X}(\widehat{\tau}^1,\widehat{\tau}^2)\ge\lambda,\,
\widehat\tau<\widehat{\sigma}_n\right)
=\mathbb{P}_{(i,\lambda)}^{\pi^1,\pi^2}
(X(\tau^1,\tau^2)\ge\lambda,\,\tau<\sigma_n).
\end{eqnarray}
The same identity holds for every strategy pair of $\mathcal G_n$ and its image under $\Psi$.
\item[(b)] If, in addition, the Lyapunov condition \eqref{eq:drift_section} holds, then for every strategy pair
$(\pi^1,\allowbreak \tau^1;\pi^2,\allowbreak \tau^2)$ of $\mathcal G$,
\begin{eqnarray}\label{bound-1}
\mathbb{P}_{(i,\lambda)}^{\pi^1,\pi^2}(\sigma_n\le\tau)\le\frac{w(i)}{n},
\end{eqnarray}
and for every strategy pair
$(\widehat\pi^1,\allowbreak \widehat\tau^1;\widehat\pi^2,\allowbreak \widehat\tau^2)$ of $\mathcal G_n$,
\begin{eqnarray}\label{bound2}
\mathbb{P}_{n,(i,\lambda)}^{\widehat{\pi}^1,\widehat{\pi}^2}(\widehat{\sigma}_n\le\widehat\tau)\le\frac{w(i)}{n}.
\end{eqnarray}
\end{itemize}
\end{proposition}
\begin{proof}
(a)
For brevity, set
$\mathbb P_n:=\allowbreak \mathbb{P}_{n,\allowbreak (i,\allowbreak \lambda)}^{\widehat{\pi}^1,\allowbreak \widehat{\pi}^2}$,
$\mathbb P:=\allowbreak \mathbb{P}_{(i,\allowbreak \lambda)}^{\pi^1,\allowbreak \pi^2}$ and
$\widehat{r}_{0:t-1}:=\allowbreak \sum\limits_{s=\allowbreak 0}^{t-1}r_n(i_s,\allowbreak a_s,\allowbreak b_s)$.

\begin{eqnarray*}
&&\mathbb P_n
\left(\widehat{X}(\widehat{\tau}^1,\widehat{\tau}^2)\ge\lambda,\,
\widehat\tau<\widehat{\sigma}_n\right)
=\sum_{t\ge0}
\mathbb P_n
\left(\widehat{X}(\widehat{\tau}^1,\widehat{\tau}^2)\ge\lambda,\,
\widehat\tau=t<\widehat{\sigma}_n\right)\\
&=&\sum_{t\ge0}
\sum_{a_0\in A_n(i_0)}\sum_{b_0\in B_n(i_0)}
\sum_{i_1\in S_n}
\cdots
\sum_{a_{t-1}\in A_n(i_{t-1})}
\sum_{b_{t-1}\in B_n(i_{t-1})}
\times\sum_{i_t\in S_n}
\prod_{s=0}^{t-1}\widehat{\pi}_s^1(a_s|\widehat h_s)
\widehat{\pi}_s^2(b_s|\widehat h_s)
p_n(i_{s+1}| i_s,a_s,b_s)\\
&&\times\Bigg[
\mathbf{1}_{\left\{\widehat{r}_{0:t-1}\ge\lambda\right\}}
\mathbb P_n(\widehat\tau_D=t\le\widehat\tau^1\wedge\widehat\tau^2
| \widehat h_t)
+\mathbf{1}_{\left\{\widehat{r}_{0:t-1}+g_n^1(i_t)\ge\lambda\right\}}
\mathbb P_n(\widehat\tau^1=t<\widehat\tau_D,\,
\widehat\tau^1\le\widehat\tau^2| \widehat h_t)\\
&&
+\mathbf{1}_{\left\{\widehat{r}_{0:t-1}+g_n^2(i_t)\ge\lambda\right\}}
\mathbb P_n(\widehat\tau^2=t<\widehat\tau_D,\,
\widehat\tau^2<\widehat\tau^1| \widehat h_t)
\Bigg]\\
&=&\mathbb P
\left(X(\tau^1,\tau^2)\ge\lambda,\,
\tau<\sigma_n\right),
\end{eqnarray*}
The last equality follows from \eqref{def:p_n},
\eqref{eq:Phi-control}, \eqref{eq:Phi-stopping},
and \eqref{eq:truncated-rewards}, since the control policies,
transition probabilities, rewards, and stopping decisions
coincide before exit.
The same argument applies to every strategy pair of
$\mathcal G_n$ and its image under $\Psi$, using
\eqref{eq:Psi-control} and \eqref{eq:Psi-stopping}.\\
	(b) Let
	$M_t:=\allowbreak w(i_t)\mathbf{1}_{\{\tau_D>t\}}$.
	By	the Lyapunov condition, for every $(i,\allowbreak \lambda)\in D^c\times  \Lambda^*$, $(\pi^1,\allowbreak \pi^2)\in \Pi^1\times \Pi^2$,
	\begin{eqnarray*}
		\mathbb{E}_{(i,\lambda)}^{\pi^1,\pi^2}[M_{t+1}|z_t,a_t,b_t]
		=\mathbb{E}_{(i,\lambda)}^{\pi^1,\pi^2}[w(i_{t+1})\mathbf{1}_{\{\tau_D>t+1\}}|z_t,a_t,b_t]
		\leq\rho w(i_t)\mathbf{1}_{\{\tau_D > t\}}\leq M_t~~~~~(\text{by}~\eqref{eq:drift_section-1}).
	\end{eqnarray*}
	Since $w\ge1$, $\{M_t\}$ is a nonnegative supermartingale.
	By Lemma~\ref{lem:ville},  we get
	\[
	\mathbb{P}_{(i,\lambda)}^{\pi^1,\pi^2}\!\left(\sup_{t\ge0}M_t\ge n\right)\le\frac{\mathbb{E}_{(i,\lambda)}^{\pi^1,\pi^2}[M_0]}{n}=\frac{w(i)}{n}.
	\]
	The event $\{\sigma_n<\tau_D\}$ implies that there exists $t<\tau_D$ with $i_t\notin S_n$, hence $w(i_t)>n$ and thus $\sup_t M_t\ge n$.
	Therefore
	$
	\{\sigma_n<\tau_D\}\subseteq\Bigl\{\sup_{t\ge0}M_t\ge n\Bigr\},\allowbreak
	$
	and consequently
	$
	\mathbb{P}_{(i,\allowbreak \lambda)}^{\pi^1,\allowbreak \pi^2}(\sigma_n<\tau_D)\le\frac{w(i)}{n}.
	$
	Since $\tau\le\tau_D$, we have
	$
	\mathbb{P}_{(i,\allowbreak \lambda)}^{\pi^1,\allowbreak \pi^2}(\sigma_n\le\tau)
	\le
	\mathbb{P}_{(i,\allowbreak \lambda)}^{\pi^1,\allowbreak \pi^2}(\sigma_n<\tau_D)
	+\mathbb{P}_{(i,\allowbreak \lambda)}^{\pi^1,\allowbreak \pi^2}(\sigma_n=\allowbreak \tau_D).
	$
	By Proposition \ref{prop:attenuation}, Assumption~\ref{ass:absorb} holds that $\mathbb{P}_{(i,\allowbreak \lambda)}^{\pi^1,\allowbreak \pi^2}(\tau_D=\allowbreak \infty)=\allowbreak 0$. Moreover, by $\{\sigma_n=\allowbreak \tau_D<\infty\}=\allowbreak \emptyset$,  we have $\mathbb{P}_{(i,\allowbreak \lambda)}^{\pi^1,\allowbreak \pi^2}(\sigma_n=\allowbreak \tau_D)=\allowbreak 0$. Thus
	$\mathbb{P}_{(i,\allowbreak \lambda)}^{\pi^1,\allowbreak \pi^2}(\sigma_n\le\tau)\le\frac{w(i)}{n}.$\\
	 Next, fix any strategy pair
$(\widehat\pi^1,\widehat\tau^1;
  \widehat\pi^2,\widehat\tau^2)$ of $\mathcal G_n$,
and let
$(\pi^1,\tau^1;\pi^2,\tau^2)
:=
\Psi(\widehat\pi^1,\widehat\tau^1;
     \widehat\pi^2,\widehat\tau^2).
$
	Since $\widehat\tau\le\widehat\tau_D$,  $\{\widehat{\sigma}_n\le\widehat\tau\}\subseteq\{\widehat{\sigma}_n\le\widehat\tau_D\}$. Thus
	\[
	\mathbb{P}_{n,(i,\lambda)}^{\widehat{\pi}^1,\widehat{\pi}^2}(\widehat{\sigma}_n\le\widehat\tau)
	\le
	\mathbb{P}_{n,(i,\lambda)}^{\widehat{\pi}^1,\widehat{\pi}^2}(\widehat{\sigma}_n<\widehat\tau_D)
	+\mathbb{P}_{n,(i,\lambda)}^{\widehat{\pi}^1,\widehat{\pi}^2}(\widehat{\sigma}_n=\widehat\tau_D).
	\]
	At any finite time, $\widehat{\sigma}_n=\allowbreak \widehat\tau_D$ is impossible, because $\widehat{\sigma}_n$ corresponds to entering the state $i^*$, whereas $i^*\notin D$.
	Hence $\mathbb{P}_{n,\allowbreak (i,\allowbreak \lambda)}^{\widehat{\pi}^1,\allowbreak \widehat{\pi}^2}(\widehat{\sigma}_n=\allowbreak \widehat\tau_D)=\allowbreak 0$, and it also holds that
	\begin{eqnarray*}
			 &&\mathbb{P}_{n,(i,\lambda)}^{\widehat{\pi}^1,\widehat{\pi}^2}(\widehat{\sigma}_n\le\widehat\tau)\\
			&\leq& \mathbb{P}_{n,(i,\lambda)}^{\widehat{\pi}^1,\widehat{\pi}^2}(\widehat{\sigma}_n<\widehat\tau_D)
			 =\sum_{t\ge1}\mathbb{P}_{n,(i,\lambda)}^{\widehat{\pi}^1,\widehat{\pi}^2}(\widehat{\sigma}_n=t,\widehat\tau_D>t)\nonumber\\	
		&=&\sum_{t\ge1}\sum_{a_0\in {A}_n(i)}\sum_{b_0\in {B}_n(i)}\sum_{i_1\in {S}_n\setminus D}\cdots\sum_{i_{t-1}\in {S}_n\setminus D}\sum_{a_{t-1}\in {A}(i_{t-1})}
		\sum_{b_{t-1}\in {B}(i_{t-1})}\widehat{\pi}_0^1(a_0|i,\lambda)\widehat{\pi}_0^2(b_0|i,\lambda){p}_n(i_1|i_0,a_0,b_0)\times\nonumber\\
		&&\cdots\times{p}_n(i_{t-1}|i_{t-2},a_{t-2},b_{t-2})
		\widehat{\pi}_{t-1}^1(a_{t-1}|{h}_{t-1})
		\times\widehat{\pi}_{t-1}^2(b_{t-1}|{h}_{t-1}){p}_n(i^*|i_{t-1},a_{t-1},b_{t-1})\nonumber\\
		&=&\sum_{t\ge1}\sum_{a_0\in A(i)}\sum_{b_0\in B(i)}\sum_{i_1\in D_n^c}\cdots\sum_{i_{t-1}\in  D_n^c}\sum_{a_{t-1}\in A(i_{t-1})}
		\sum_{b_{t-1}\in B(i_{t-1})}\sum_{i_t\notin {S}_n}{\pi}_0^1(a_0|i,\lambda){\pi}_0^2(b_0|i,\lambda)p(i_1|i_0,a_0,b_0)\times\nonumber\\
		&&\cdots\times p(i_{t-1}|i_{t-2},a_{t-2},b_{t-2})
		{\pi}_{t-1}^1(a_{t-1}|h_{t-1})\times {\pi}_{t-1}^2(b_{t-1}|{h}_{t-1}){p}(i_t|i_{t-1},a_{t-1},b_{t-1})\qquad\text{(by  \eqref{def:p_n})}\nonumber\\
		&=&\sum_{t\ge1}\mathbb{P}_{(i,\lambda)}^{{\pi}^1,{\pi}^2}({\sigma}_n=t,{\tau}_D>t)
		=\mathbb{P}_{(i,\lambda)}^{\pi^1,\pi^2}(\sigma_n<\tau_D)\le\frac{w(i)}{n},\nonumber
	\end{eqnarray*}
	which completes Proposition \ref{prop:two-estimates}.
\end{proof}

We combine the two estimates  on $\{\widehat\tau<\widehat\sigma_n\}$ and $\{\widehat\sigma_n\le\widehat\tau\}$. Integrating the pointwise bounds with respect to $\nu$ yields error bounds between the values of $\mathcal G_n$ and $\mathcal G$, and shows that approximate Nash equilibria of the truncated games remain approximate Nash equilibria of the original game under the map $\Psi$.
\begin{theorem}\label{thm:truncation}
	Suppose the Lyapunov condition \eqref{eq:drift_section} holds.
	Let $\nu\in\mathcal{P}(S)$ be such that
	$
	\nu(w):=\allowbreak \sum_{i\in D^c}w(i)\nu(i)<\infty.$ Then the following statements hold.
	\begin{itemize}
		\item[(a)] The truncated game values satisfy
			\begin{eqnarray}
				\bigl|J_{n,\nu_n}^*-J_\nu^*\bigr|\le \frac{9\nu(w)}{n} \ \ \ \ \forall n\geq 1.
			\label{eq:truncation-value-bound}
			\end{eqnarray}
			Consequently, $\lim\limits_{n\to\infty}J_{n,\allowbreak \nu_n}^*=\allowbreak J_\nu^*$.
			\item[(b)] If  $(\widehat{\pi}^1_\epsilon,\allowbreak \widehat{\tau}^1_\epsilon;\widehat{\pi}^2_\epsilon,\allowbreak \widehat{\tau}^2_\epsilon)$  is an $\epsilon$-Nash equilibrium of $\mathcal{G}_n$ with respect to $(\nu_n,\allowbreak \lambda^*)$, then $\Psi(\widehat{\pi}^1_\epsilon,\allowbreak \widehat{\tau}^1_\epsilon;\widehat{\pi}^2_\epsilon,\allowbreak \widehat{\tau}^2_\epsilon)$ is an $(\epsilon+18\nu(w)/n)$-Nash equilibrium of $\mathcal{G}$ with respect to $(\nu,\allowbreak \lambda^*)$.
	\end{itemize}
\end{theorem}
\begin{proof}  For any $(i,\allowbreak \lambda)\in  D_n^c\times\Lambda_n^*$ and  $(\pi^1,\allowbreak \tau^1,\allowbreak \pi^2,\allowbreak \tau^2)\in (\Pi^1\times\Gamma)\times(\Pi^2\times\Gamma)$, let $(\widehat{\pi}^1,\allowbreak \widehat{\tau}^1;\widehat{\pi}^2,\allowbreak \widehat{\tau}^2):=\allowbreak \Phi(\pi^1,\allowbreak \tau^1;\pi^2,\allowbreak \tau^2)$. Then, by Proposition~\ref{prop:two-estimates}, we have
			\begin{eqnarray}
				&&\bigl|\mathbb{P}_{(i,\lambda)}^{\pi^1,\pi^2}(X(\tau^1,\tau^2)\ge\lambda) - \mathbb{P}_{n,(i,\lambda)}^{\widehat{\pi}^1,\widehat{\pi}^2}(\widehat{X}(\widehat{\tau}^1,\widehat{\tau}^2)\ge\lambda)\bigr|\nonumber\\
				 &\le&\Big|\mathbb{P}_{(i,\lambda)}^{\pi^1,\pi^2}(X(\tau^1,\tau^2)\ge\lambda,\tau<\sigma_n)-\mathbb{P}_{n,(i,\lambda)}^{\widehat{\pi}^1,\widehat{\pi}^2}(\widehat{X}(\widehat{\tau}^1,\widehat{\tau}^2)\ge\lambda,
			\widehat\tau<\widehat{\sigma}_n)\Big|+\mathbb{P}_{(i,\lambda)}^{\pi^1,\pi^2}(\sigma_n\le\tau)	 +\mathbb{P}_{n,(i,\lambda)}^{\widehat{\pi}^1,\widehat{\pi}^2}(\widehat{\sigma}_n\le\widehat\tau)\nonumber\\
			&	=&\mathbb{P}_{(i,\lambda)}^{\pi^1,\pi^2}(\sigma_n\le\tau)
				 +\mathbb{P}_{n,(i,\lambda)}^{\widehat{\pi}^1,\widehat{\pi}^2}(\widehat{\sigma}_n\le\widehat\tau)~~~(\text{by}~\eqref{eq:agg-noexit})\label{coupling1-1}\\
				 &\le&\frac{2w(i)}{n}.~~~~~~~~~~~~~~~~~~~~~~~~~~(\text{by}~\eqref{bound-1}-\eqref{bound2})\label{coupling2}
			\end{eqnarray}
			Similarly, let $(\pi^1,\allowbreak \tau^1;\pi^2,\allowbreak \tau^2):=\allowbreak \Psi(\widehat{\pi}^1,\allowbreak \widehat{\tau}^1;\widehat{\pi}^2,\allowbreak \widehat{\tau}^2)$ for any  $(\widehat{\pi}^1,\allowbreak \widehat{\tau}^1;\widehat{\pi}^2,\allowbreak \widehat{\tau}^2) \in (\Pi^1_n\times\Gamma_n)\times(\Pi^2_n\times\Gamma_n)$. Then,
		\begin{eqnarray*}
			 \bigl|\mathbb{P}_{n,(i,\lambda)}^{\widehat{\pi}^1,\widehat{\pi}^2}(\widehat{X}(\widehat{\tau}^1,\widehat{\tau}^2)\ge\lambda)-\mathbb{P}_{(i,\lambda)}^{\pi^1,\pi^2}(X(\tau^1,\tau^2)\ge\lambda)\bigr|\le\frac{2w(i)}{n}.
		\end{eqnarray*}
		Indeed, write
		$
		V_i:=\allowbreak \mathbb{P}_{(i,\allowbreak \lambda^*)}^{\pi^1,\allowbreak \pi^2}(X(\tau^1,\allowbreak \tau^2)\ge\lambda^*)
		$ for $i\in D^c$
		and
		$
		\widehat V_i:=\allowbreak \mathbb{P}_{n,\allowbreak (i,\allowbreak \lambda^*)}^{\widehat{\pi}^1,\allowbreak \widehat{\pi}^2}
		(\widehat{X}(\widehat{\tau}^1,\allowbreak \widehat{\tau}^2)\ge\lambda^*)
		$ for $i\in  D_n^c$.
		Since $\nu_n(i)=\allowbreak \nu(i)$ for $i\in  D_n^c$ and the mass $\nu(D^c\setminus S_n)$ is sent to state $i^*$, whose payoff is zero when $\lambda^*>0$, we have
		\begin{eqnarray}\label{value-bound}
&&\bigl|J_{\nu}(\pi^1,\tau^1;\pi^2,\tau^2)
		-J_{n,\nu_n}(\widehat{\pi}^1,\widehat{\tau}^1;\widehat{\pi}^2,\widehat{\tau}^2)\bigr|
=\left|\sum_{i\in  D_n^c}(V_i-\widehat V_i)\nu(i)
		+\sum_{i\in D^c\setminus S_n}V_i\nu(i)\right|\nonumber\\
	&	\le&\sum_{i\in  D_n^c}|V_i-\widehat V_i|\nu(i)
		+\sum_{i\in D^c\setminus S_n}V_i\nu(i)\nonumber\\
	&	\le&\sum_{i\in  D_n^c}\frac{2w(i)}{n}\nu(i)+\nu(D^c\setminus S_n).~(\text{by}~\eqref{coupling2} ~\text{and}~ 0\le V_i\le1)~~~~
		\end{eqnarray}
	Moreover,by $\sum\limits_{i\in  D_n^c}\frac{2w(i)}{n}\nu(i)\le\frac{2\nu(w)}{n},\allowbreak $  $D^c\setminus S_n\subset\{i\in D^c:w(i)>n\}$, and
		$\mathbf{1}_{\{w(i)>n\}}\le w(i)/n$, we have
		\[
		\nu(D^c\setminus S_n)
		\le \sum_{i\in D^c}\mathbf{1}_{\{w(i)>n\}}\nu(i)
		\le \frac{1}{n}\sum_{i\in D^c}w(i)\nu(i)
		=\frac{\nu(w)}{n}.
		\]
		Therefore, 	by \eqref{value-bound}, we have
		\begin{eqnarray}\label{dist-coupling-phi}
\bigl|J_{\nu}(\pi^1,\tau^1;\pi^2,\tau^2)
		-J_{n,\nu_n}(\widehat{\pi}^1,\widehat{\tau}^1; \widehat{\pi}^2,\widehat{\tau}^2)\bigr|
		\le \frac{3\nu(w)}{n}.
		\end{eqnarray}
		Similarly, for  any  $(\widehat{\pi}^1,\allowbreak \widehat{\tau}^1;\widehat{\pi}^2,\allowbreak \widehat{\tau}^2) \in (\Pi^1_n\times\Gamma_n)\times(\Pi^2_n\times\Gamma_n)$, $(\pi^1,\allowbreak \tau^1;\pi^2,\allowbreak \tau^2) =\allowbreak  \Psi(\widehat{\pi}^1,\allowbreak \widehat{\tau}^1;\widehat{\pi}^2,\allowbreak \widehat{\tau}^2)$, we have
		\begin{eqnarray}\label{dist-coupling-psi}
\bigl|J_{n,\nu_n}(\widehat{\pi}^1,\widehat{\tau}^1;\widehat{\pi}^2,\widehat{\tau}^2)
		-J_{\nu}(\pi^1,\tau^1;\pi^2,\tau^2)\bigr|
		\le \frac{3\nu(w)}{n}.
		\end{eqnarray}
	Using the obtained facts above, we next prove (a) and (b).

		(a) Fix $\epsilon>0$, and let
	$(\widehat{\pi}^1_\epsilon,\allowbreak \widehat{\tau}^1_\epsilon;\widehat{\pi}^2_\epsilon,\allowbreak \widehat{\tau}^2_\epsilon)$
	be an $\epsilon$-Nash equilibrium of $\mathcal{G}_n$ with respect to $(\nu_n,\allowbreak \lambda^*)$. Set
	$
	(\pi^1_\epsilon,\allowbreak \tau^1_\epsilon;\pi^2_\epsilon,\allowbreak \tau^2_\epsilon)
	:=\allowbreak \Psi(\widehat{\pi}^1_\epsilon,\allowbreak \widehat{\tau}^1_\epsilon;\widehat{\pi}^2_\epsilon,\allowbreak \widehat{\tau}^2_\epsilon).
	$
Moreover, for any $(\pi^1,\allowbreak \tau^1)\in\Pi^1\times\Gamma$, applying \eqref{dist-coupling-phi} and \eqref{dist-coupling-psi}, we obtain
		\begin{eqnarray}
			&&J_\nu(\pi^1,\tau^1;\pi^2_\epsilon,\tau^2_\epsilon)
		=J_\nu(\pi^1,\tau^1;\Psi_2(\widehat{\pi}^2_\epsilon,\widehat{\tau}^2_\epsilon))\le J_{n,\nu_n}(\Phi_1({\pi}^1,{\tau}^1);\Phi_2(\Psi_2(\widehat{\pi}^2_\epsilon, \widehat{\tau}^2_\epsilon)))+\frac{3\nu(w)}{n}\nonumber\\
		 &\le&J_\nu\big(\Psi_1(\Phi_1({\pi}^1,{\tau}^1));\Psi_2(\Phi_2(\Psi_2(\widehat{\pi}^2_\epsilon,\widehat{\tau}^2_\epsilon)))\big)
		+\frac{6\nu(w)}{n}\nonumber\\
		 &=&J_\nu(\Psi_1(\Phi_1({\pi}^1,{\tau}^1));\Psi_2(\widehat{\pi}^2_\epsilon,\widehat{\tau}^2_\epsilon))+\frac{6\nu(w)}{n}\quad(\text{by \eqref{eq:PhiPsi}})\\
		& \le&J_{n,\nu_n}(\Phi_1({\pi}^1,{\tau}^1);(\widehat{\pi}^{2}_\epsilon,\widehat{\tau}^{2}_\epsilon))+\frac{9\nu(w)}{n}\label{eq:truncation-psi-step}\\
	&\le&J_{n,\nu_n}^*+\epsilon+\frac{9\nu(w)}{n},\label{eq:truncation-epsilon-Nash}
		\end{eqnarray}
		where \eqref{eq:truncation-psi-step} is due to \eqref{dist-coupling-psi} that $J_{n,\allowbreak \nu_n}(\Phi_1({\pi}^1,\allowbreak {\tau}^1);(\widehat{\pi}^{2}_\epsilon,\allowbreak \widehat{\tau}^{2}_\epsilon))\ge J_\nu(\Psi_1(\Phi_1({\pi}^1,\allowbreak {\tau}^1));\Psi_2(\widehat{\pi}^2_\epsilon,\allowbreak \widehat{\tau}^2_\epsilon))-\frac{3\nu(w)}{n}$, and \eqref{eq:truncation-epsilon-Nash} follows from the $\epsilon$-optimality of $(\widehat{\pi}^{2}_\epsilon,\allowbreak \widehat{\tau}^{2}_\epsilon)$ in $\mathcal{G}_n$.\\
		Taking the supremum over $(\pi^1,\allowbreak \tau^1)$ on the left-hand side gives
		\begin{eqnarray}\label{upper-deviation}
		\sup_{(\pi^1,\tau^1)}J_\nu(\pi^1,\tau^1;\pi^2_\epsilon,\tau^2_\epsilon)
		\le J_{n,\nu_n}^*
		+\epsilon+\frac{9\nu(w)}{n},
		\end{eqnarray}
		which implies that
		\begin{eqnarray}\label{upper-deviation-1}
			\inf_{(\pi^2,\tau^2)}\sup_{(\pi^1,\tau^1)}
			J_\nu(\pi^1,\tau^1;\pi^2,\tau^2)
			\le J_{n,\nu_n}^*
		+\epsilon+\frac{9\nu(w)}{n},
		\end{eqnarray}
		Similarly, we can show the reverse inequality:
		\begin{eqnarray}\label{lower-deviation}	
			\inf_{(\pi^2,\tau^2)}J_\nu(\pi^1_\epsilon,\tau^1_\epsilon;\pi^2,\tau^2)
		\ge J_{n,\nu_n}^*-\epsilon-\frac{9\nu(w)}{n},
		\end{eqnarray}
	which implies that
			\begin{eqnarray}\label{lower-deviation-1}
		\sup_{(\pi^1,\tau^1)}\inf_{(\pi^2,\tau^2)}
		J_\nu(\pi^1,\tau^1;\pi^2,\tau^2)
		\ge J_{n,\nu_n}^*
		-\epsilon-\frac{9\nu(w)}{n}.
		\end{eqnarray}
		By \eqref{lower-deviation-1} and \eqref{upper-deviation-1}, we have
		$
		|J_\nu^*-J_{n,\allowbreak \nu_n}^*|\le\epsilon+\frac{9\nu(w)}{n}.
		$
		Letting $\epsilon\to0$ yields $|J_\nu^*-J_{n,\allowbreak \nu_n}^*|\le\frac{9\nu(w)}{n}$, which implies $\lim\limits_{n\to\infty}J_{n,\allowbreak \nu_n}^*=\allowbreak J_\nu^*$, i.e.,  (a) holds.

(b) Using (a) together with \eqref{lower-deviation} and \eqref{upper-deviation}, we conclude that $(\pi^1_\epsilon,\allowbreak \tau^1_\epsilon;\pi^2_\epsilon,\allowbreak \tau^2_\epsilon)$ is an $(\epsilon+18\nu(w)/n)$-Nash equilibrium of $\mathcal{G}$ with respect to $(\nu,\allowbreak \lambda^*)$, which proves (b).
	\end{proof}
\begin{remark}\label{rem:truncation-novelty}
The finite-approximation results in
\citet{DufourPrietoRumeau2019,TidballAltman1996,
LorenzoHernandezPrietoRumeau2015,TidballPourtallierAltman1997}
concern risk-neutral zero-sum games.
Theorem~\ref{thm:truncation} here establishes explicit
 error bounds between the values of the  truncated and original games for zero-sum
stopping games under the probability criterion.
To the best of our knowledge, these error bounds are new.
\end{remark}

\section[An algorithm for epsilon-Nash Equilibria]{An Algorithm for $\epsilon$-Nash Equilibria}\label{sec:algorithm}
In this section, we present an algorithm for computing $\epsilon$-Nash equilibria of the original game $\mathcal{G}$. The algorithm is based on the error bounds in Theorem~\ref{thm:truncation}, and works under the Lyapunov condition \eqref{eq:drift_section} with $
	\nu(w)<\infty.$ For finite-dimensional implementation, we further assume that the Lyapunov sublevel sets  in \eqref{sublevel-set} are finite, i.e.,
\begin{eqnarray}\label{eq:finite-sublevel}
D_n^c=\{\,i\in D^c:w(i)\le n\,\}\quad\text{is finite for every }n\in\mathbb{N}.
\end{eqnarray}
If $D_n^c=\varnothing$, computation is trivial.
Otherwise, finiteness of the admissible action sets and
positivity of $r_n$ on $D_n^c$ imply that the set
$
\mathcal R_n
:=
\{r_n(i,a,b):i\in D_n^c,\ a\in A_n(i),\ b\in B_n(i)\}
$
is finite and has a positive minimum
$r_{n,\allowbreak \min}:=\allowbreak \min\mathcal R_n>0$.
Since $\mathcal T_nu(i,\allowbreak \lambda)=\allowbreak 1$ for $\lambda<0$, it suffices to
consider
$
\bar\Lambda_n^*:=\Lambda_n^*\cap[0,\lambda^*].
$
Every threshold in $\bar\Lambda_n^*$ is obtained by subtracting at most
$\lfloor\lambda^*/r_{n,\allowbreak \min}\rfloor$ elements of the finite set
$\mathcal R_n$ from $\lambda^*$. Hence $\bar\Lambda_n^*$ is finite, and
so is
$
E_n:=\allowbreak D_n^c\times\bar\Lambda_n^*.$   Hence, under
the finite-sublevel condition, each truncated game has an equivalent
finite-state representation.
Then, we state our algorithm (say, Algorithm~\ref{alg:truncVI}) below.
\\
\begin{breakablealgorithm}
\caption{Computation of an $\epsilon$-Nash Equilibrium for the Original Game}
\label{alg:truncVI}
\begin{algorithmic}
\State \textbf{Step 1 (Input and truncation level).}
Given initial distribution $\nu\in\mathcal{P}(S)$ with $\nu(w)<\infty$, threshold $\lambda^*>0$, and accuracy $\epsilon>0$.
Choose an integer
$
n\ge\left\lceil\frac{24\nu(w)}{\epsilon}\right\rceil.
$
Construct the truncated state set $S_n:=\allowbreak \{i\in S:w(i)\le n\}\cup D$, and define the truncated initial distribution $\nu_n$ by $\nu_n(i):=\allowbreak \nu(i)$ for $i\in S_n$ and $\nu_n(i^*):=\allowbreak \nu(S\setminus S_n)$.

\State \textbf{Step 2 (Initialization).}
Set $u_{n,\allowbreak 0}^0(i,\allowbreak \lambda)=\allowbreak \mathbf{1}_{\{\lambda\le0\}}$ and $u_{n,\allowbreak 0}^1(i,\allowbreak \lambda)=\allowbreak 1$ for all $(i,\allowbreak \lambda)\in E_n$, $u_{n,\allowbreak t}^l(i,\allowbreak \lambda)=\allowbreak \mathbf{1}_{\{\lambda\le0\}}$ for all $(i,\allowbreak \lambda)\in D^*\times\bar{\Lambda}_n^*$, $u_{n,\allowbreak t}^l(i,\allowbreak \lambda)=\allowbreak 1$ for all $i\in D_n^c,\allowbreak  \lambda<0$, $t\ge0$, $l=\allowbreak 0,\allowbreak 1$.
Let $m=\allowbreak 0$.

\State \textbf{Step 3 (Value iteration).}
\Repeat
    \For{each $(i,\allowbreak \lambda)\in E_n$ and $l=\allowbreak 0,\allowbreak 1$}
        \State Compute the payoff matrix with entries
        $\mathbb{G}_n u_{n,\allowbreak m}^l(i,\allowbreak \lambda;a,\allowbreak b)$ for all $(a,\allowbreak b)\in A_n(i)\times B_n(i)$.
        \State Solve the zero-sum matrix game with these payoffs to obtain its value $\overline{u}_{n,\allowbreak m}^l(i,\allowbreak \lambda)$.
        \State Update
        \[
        u_{n,m+1}^l(i,\lambda)
        :=
        \mathbf{1}_{\{g^1_n(i)\ge\lambda\}}\vee
        \bigl(\mathbf{1}_{\{g^2_n(i)\ge\lambda\}}
        \wedge\overline{u}_{n,m}^l(i,\lambda)\bigr).
        \]
    \EndFor
    \State $m\gets m+1$
\Until{$u_{n,\allowbreak m}^1-u_{n,\allowbreak m}^0\le\epsilon/4$ on $E_n$} (on-line stopping rule)
\State Set $m_\epsilon:=\allowbreak m$
\begin{eqnarray}\label{m_epsilon}
(i.e., ~m_\epsilon:=\min\{m: u_{n,m}^1-u_{n,m}^0\le\epsilon/4\}~\text{on}~E_n).
\end{eqnarray}
\State \textbf{Step 4 (An $\epsilon/4$-Nash equilibrium of the truncated game $\mathcal{G}_n$).}
\For{each $(i,\allowbreak \lambda)\in E_n$}
    \State Solve the zero-sum matrix game with payoffs
    $\mathbb{G}_n u_{n,\allowbreak m_\epsilon}^0(i,\allowbreak \lambda;\cdot,\allowbreak \cdot)$, and let
    $\widehat{\pi}_{n,\allowbreak m_\epsilon}^1(\cdot|i,\allowbreak \lambda)$ be an optimal mixed action for Player~1.
    \State Solve the zero-sum matrix game with payoffs
    $\mathbb{G}_n u_{n,\allowbreak m_\epsilon}^1(i,\allowbreak \lambda;\cdot,\allowbreak \cdot)$, and let
    $\widehat{\pi}_{n,\allowbreak m_\epsilon}^2(\cdot|i,\allowbreak \lambda)$ be an optimal mixed action for Player~2.
 \EndFor
	\State Define stopping times
    \begin{eqnarray}\label{eq:truncation-stopping-times}
    \mspace{-20mu}\widehat{\tau}_{n,m_\epsilon}^1
    :=\inf\Bigl\{t\ge0:
    \mathcal{L}_n u_{n,m_\epsilon}^0(\widehat i_t,\widehat\lambda_t)
    \le
    \mathbf{1}_{\{g^1_n(\widehat i_t)\ge\widehat\lambda_t\}}\Bigr\},
	\end{eqnarray}
    \begin{eqnarray}\label{eq:truncation-stopping-times-2}
     \mspace{-20mu}\widehat{\tau}_{n,m_\epsilon}^2
    :=\inf\Bigl\{t\ge0:
    \mathcal{L}_n u_{n,m_\epsilon}^1(\widehat i_t,\widehat\lambda_t)
    \ge
    \mathbf{1}_{\{g^2_n(\widehat i_t)\ge\widehat\lambda_t\}}\Bigr\}.
	\end{eqnarray}
\State Output the $\epsilon/4$-Nash equilibrium of the truncated game $\mathcal{G}_n$
\begin{eqnarray}\label{epsilon-Nash-truncated}
\bigl(
\widehat{\pi}_{n,m_\epsilon}^1,\widehat{\tau}_{n,m_\epsilon}^1;
\widehat{\pi}_{n,m_\epsilon}^2,\widehat{\tau}_{n,m_\epsilon}^2
\bigr).
\end{eqnarray}

\State \textbf{Step 5  (An $\epsilon$-Nash equilibrium of the original game $\mathcal{G}$).}\\
 Map back and output the $\epsilon$-Nash equilibrium  and the value  of the original game $\mathcal{G}$:
\begin{align}
(\pi_{n,m_\epsilon}^1,\tau_{n,m_\epsilon}^1;
\pi_{n,m_\epsilon}^2,\tau_{n,m_\epsilon}^2)&:=\Psi(\widehat{\pi}_{n,m_\epsilon}^1,\widehat{\tau}_{n,m_\epsilon}^1;
\widehat{\pi}_{n,m_\epsilon}^2,\widehat{\tau}_{n,m_\epsilon}^2
),\label{epsilon-Nash-orignal}\\
J_{n,\nu_n,m_\epsilon}
&:=\sum_{i\in D_n^c} u_{n,m_\epsilon}^0(i,\lambda^*)\nu_n(i).\label{epsilon-value-orignal}
\end{align}
\end{algorithmic}
\end{breakablealgorithm}
\smallskip
\begin{remark} Given any truncated level $n\geq 1$, the stopping rule of Algorithm~\ref{alg:truncVI} is based on the difference $u_{n,\allowbreak m}^1-u_{n,\allowbreak m}^0$ of the two sequences $u_{n,\allowbreak m}^1$ and $u_{n,\allowbreak m}^0$, which does not require the contraction
condition~\eqref{one-step-absorb} used in
\citet{Huang2011,HG-2020,Liu2001}.
\end{remark}
The following theorem establishes the finite termination and convergence of Algorithm~\ref{alg:truncVI} on $\mathcal G_n$ as $\epsilon\to0$.
\begin{theorem}\label{thm:truncated-convergence} Suppose that \eqref{eq:drift_section} and \eqref{eq:finite-sublevel}
hold.
Fix $n\in\mathbb N$. For each $\epsilon>0$, let $m_\epsilon$ be the
integer determined by \eqref{m_epsilon} after running Steps 2--4 of Algorithm~\ref{alg:truncVI} on
$\mathcal G_n$. Then, we have the following assertions:
\begin{itemize}
\item[(a)]
 The output $(\widehat{\pi}_{n,\allowbreak m_\epsilon}^1,\allowbreak \widehat{\tau}_{n,\allowbreak m_\epsilon}^1;
\widehat{\pi}_{n,\allowbreak m_\epsilon}^2,\allowbreak \widehat{\tau}_{n,\allowbreak m_\epsilon}^2
)$ in
\eqref{epsilon-Nash-truncated} is obtained  after finitely many iterations of Algorithm 1, and is an $\epsilon/4$-Nash equilibrium of
$\mathcal G_n$ with respect to $(\nu_n,\allowbreak \lambda^*)$;
\item[(b)] Let $\{\epsilon_k\}_{k\ge1}$ satisfy
$\epsilon_k\downarrow0$, and set $m_k:=\allowbreak m_{\epsilon_k}$. Then
\begin{itemize}
\item[(i)]the policy sequence $\{(\widehat{\pi}_{n,\allowbreak m_k}^1,\allowbreak \widehat{\pi}_{n,\allowbreak m_k}^2)\}_{k\ge1}$
has an accumulation point denoted by $(\widehat{\pi}_n^{*,\allowbreak 1},\allowbreak \widehat{\pi}_n^{*,\allowbreak 2})$;
\item[(ii)]the stopping-time sequence $\{(\widehat{\tau}_{n,\allowbreak m_k}^1,\allowbreak \widehat{\tau}_{n,\allowbreak m_k}^2)\}_{k\ge1}$
converges to $(\widehat{\tau}_{u_n^*}^1,\allowbreak \widehat{\tau}_{u_n^*}^2)$ defined in \eqref{eq:truncation-stopping-times}--\eqref{eq:truncation-stopping-times-2}.
\end{itemize}
Consequently,
$
\bigl(
\widehat{\pi}_n^{*,\allowbreak 1},\allowbreak \widehat{\tau}_{u_n^*}^1;
\widehat{\pi}_n^{*,\allowbreak 2},\allowbreak \widehat{\tau}_{u_n^*}^2
\bigr)
$
is a Nash equilibrium of $\mathcal G_n$ with respect to
$(\nu_n,\allowbreak \lambda^*)$.
\end{itemize}
\end{theorem}
\begin{proof}
(a) By Corollary~\ref{cor:truncated-game-Nash},
for every $m\ge1$, we have
$u_{n,\allowbreak m}^0\le u_n^*\le u_{n,\allowbreak m}^1$ on $E_n,\allowbreak
$
and $u_n^*=\allowbreak \mathcal L_nu_n^*$. Applying Theorem~\ref{thm:relationship} to $\mathcal G_n$, for the output $
\bigl(
\widehat{\pi}_{n,\allowbreak m_\epsilon}^1,\allowbreak \widehat{\tau}_{n,\allowbreak m_\epsilon}^1;
\widehat{\pi}_{n,\allowbreak m_\epsilon}^2,\allowbreak \widehat{\tau}_{n,\allowbreak m_\epsilon}^2
\bigr)$ in
\eqref{epsilon-Nash-truncated}, we obtain that
\begin{eqnarray}\label{convergence-ineq1}
J_{n,\nu_n}(\widehat{\pi}_{n,m_\epsilon}^1,\widehat{\tau}_{n,m_\epsilon}^1;\widehat{\pi}^2,\widehat{\tau}^2)
\ge
\sum_{i\in D_n^c}u_{n,m_\epsilon}^0(i,\lambda^*)\nu_n(i)~~\forall (\widehat{\pi}^2,\widehat{\tau}^2)\in \Pi^2_n\times \Gamma_n,
\end{eqnarray}
\begin{eqnarray}\label{convergence-ineq2}
J_{n,\nu_n}(\widehat{\pi}^1,\widehat{\tau}^1;\widehat{\pi}_{n,m_\epsilon}^2,\widehat{\tau}_{n,m_\epsilon}^2)
\le
\sum_{i\in D_n^c}u_{n,m_\epsilon}^1(i,\lambda^*)\nu_n(i)~~\forall (\widehat{\pi}^1,\widehat{\tau}^1)\in \Pi^1_n\times \Gamma_n.
\end{eqnarray}
By the stopping rule in Step 3 of Algorithm~\ref{alg:truncVI}, we have
$
u_{n,\allowbreak m_\epsilon}^1(i,\allowbreak \lambda)-u_{n,\allowbreak m_\epsilon}^0(i,\allowbreak \lambda)\le\frac{\epsilon}{4}
$ for all $(i,\allowbreak \lambda)\in E_n$, together with $u_{n,\allowbreak m_\epsilon}^0\le u_n^*\le u_{n,\allowbreak m_\epsilon}^1$, which yields
\begin{eqnarray*}
\sum_{i\in D_n^c}u_{n,m_\epsilon}^0(i,\lambda^*)\nu_n(i)
\ge J_{n,\nu_n}^*-\frac{\epsilon}{4}
\qquad\text{and}\qquad
\sum_{i\in D_n^c}u_{n,m_\epsilon}^1(i,\lambda^*)\nu_n(i)
\le J_{n,\nu_n}^*+\frac{\epsilon}{4}.
\end{eqnarray*}
Combining with \eqref{convergence-ineq1} and \eqref{convergence-ineq2}, we obtain
$\bigl(
\widehat{\pi}_{n,\allowbreak m_\epsilon}^1,\allowbreak \widehat{\tau}_{n,\allowbreak m_\epsilon}^1;
\widehat{\pi}_{n,\allowbreak m_\epsilon}^2,\allowbreak \widehat{\tau}_{n,\allowbreak m_\epsilon}^2
\bigr)$ is an $\epsilon/4$-Nash equilibrium of $\mathcal G_n$.

(b)
(i) Since $E_n$, $A_n(i)$ and $B_n(i)$ are finite,
$
\mathscr K_n
:=\allowbreak
\prod_{(i,\allowbreak \lambda)\in E_n}\mathcal P(A_n(i))
\times
\prod_{(i,\allowbreak \lambda)\in E_n}\mathcal P(B_n(i))
$
is compact. Hence the policy sequence$
\bigl\{(\widehat{\pi}_{n,\allowbreak m_k}^1,\allowbreak \widehat{\pi}_{n,\allowbreak m_k}^2)\bigr\}_{k\ge1}
\subset\mathscr K_n$
has an accumulation point, denoted by
$(\widehat{\pi}_n^{*,\allowbreak 1},\allowbreak \widehat{\pi}_n^{*,\allowbreak 2})$.
Passing to a subsequence and keeping the same index notation, we may assume
\[
(\widehat{\pi}_{n,m_k}^1,\widehat{\pi}_{n,m_k}^2)
\longrightarrow
(\widehat{\pi}_n^{*,1},\widehat{\pi}_n^{*,2})
\qquad\text{as }k\to\infty.
\]
For each $k$, Step~4 of Algorithm~\ref{alg:truncVI} gives
\begin{eqnarray*}
\mathcal T_nu_{n,m_k}^0(i,\lambda)
=
\inf_{\psi\in\mathcal P(B_n(i))}
\sum_{a,b}\mathbb G_nu_{n,m_k}^0(i,\lambda;a,b)
\widehat{\pi}_{n,m_k}^1(a|i,\lambda)\psi(b).
\end{eqnarray*}
Letting $k\to\infty$ and using $u_{n,\allowbreak m_k}^0\to u_n^*$ pointwise yields
\[
\mathcal T_nu_n^*(i,\lambda)
=
\inf_{\psi\in\mathcal P(B_n(i))}
\sum_{a,b}\mathbb G_nu_n^*(i,\lambda;a,b)
\widehat{\pi}_n^{*,1}(a|i,\lambda)\psi(b).
\]
Similarly,
\[
\mathcal T_nu_n^*(i,\lambda)
=
\sup_{\phi\in\mathcal P(A_n(i))}
\sum_{a,b}\mathbb G_nu_n^*(i,\lambda;a,b)
\phi(a)\widehat{\pi}_n^{*,2}(b|i,\lambda).
\]
Thus $(\widehat{\pi}_n^{*,\allowbreak 1},\allowbreak \widehat{\pi}_n^{*,\allowbreak 2})$ is a saddle point of $\mathbb G_nu_n^*$.
\\
(ii) Let $m_k:=\allowbreak m_{\epsilon_k}$. Since $m_\epsilon$ is the smallest integer satisfying
$0\leq u_{n,\allowbreak m}^1-u_{n,\allowbreak m}^0\leq \epsilon/4$ on $E_n$, and
$\epsilon_k\downarrow0$, the sequence $\{m_k\}$ is non-decreasing.
If it is bounded, then $m_k=\allowbreak \bar m$ eventually, and its definition
implies
$u_{n,\allowbreak \bar m}^0=\allowbreak u_{n,\allowbreak \bar m}^1=\allowbreak u_n^*$.
It therefore remains to consider the case of $m_k\to\infty$.

By \eqref{eq:truncation-stopping-times} and the monotone convergence
$\mathcal L_n u_{n,\allowbreak m_k}^0\uparrow u_n^*$, the
following set for Player~1,
\[
\bigl\{(i,\lambda)\in E_n:
\mathcal L_n u_{n,m_k}^0(i,\lambda)\le
\mathbf 1_{\{g_n^1(i)\ge\lambda\}}\bigr\},
\]
is non-increasing in $m_k$. Hence $\widehat\tau_{n,\allowbreak m_k}^1$ is non-decreasing in $k$ and converges
pointwise to $\widehat\tau_{u_n^*}^1$ as $k\to\infty$.
Similarly, since $\mathcal L_n u_{n,\allowbreak m_k}^1\downarrow u_n^*$, we see that
$\widehat\tau_{n,\allowbreak m_k}^2$  is non-decreasing in $m_k$ and converges pointwise to
$\widehat\tau_{u_n^*}^2$. Thus,   $\{\widehat\tau_{n,\allowbreak m_k}^1,\allowbreak \widehat\tau_{n,\allowbreak m_k}^2\}_{k\ge1}$ converges pointwise to
$(\widehat{\tau}_{u_n^*}^1,\allowbreak \widehat{\tau}_{u_n^*}^2)$, where $\widehat\tau_{u_n^*}^1$ and $\widehat\tau_{u_n^*}^2$
are the stopping times induced by $u_n^*$ through \eqref{eq:truncation-stopping-times}-\eqref{eq:truncation-stopping-times-2}.

Combined with (i) and (ii), by Corollary~\ref{cor:truncated-game-Nash}(c), we conclude that $(\widehat{\pi}_n^{*,\allowbreak 1},\allowbreak \widehat{\tau}_{u_n^*}^1;
\widehat{\pi}_n^{*,\allowbreak 2},\allowbreak \widehat{\tau}_{u_n^*}^2
\bigr)$ (induced by $u_n^*$) is a Nash
equilibrium of $\mathcal G_n$ with respect to $(\nu_n,\allowbreak \lambda^*)$.
\end{proof}

We next study the convergence of the $\epsilon_k$-Nash
equilibria obtained by Algorithm~\ref{alg:truncVI}
as the truncation level tends to infinity.
The topology on the strategy space used in the following
theorem will be introduced below.

\begin{theorem}\label{thm:main_epsilon}
	Suppose that the Lyapunov condition \eqref{eq:drift_section} and
the finite-sublevel condition \eqref{eq:finite-sublevel} hold. Let
$\nu\in\mathcal P(S)$ satisfy $\nu(w)<\infty$, and let
$\lambda^*>0$ and $\epsilon>0$. Choose $n\in\mathbb N$ such that
$
n\ge
\left\lceil\frac{24\nu(w)}{\epsilon}\right\rceil,
$
and run Algorithm~\ref{alg:truncVI}. Then:
\begin{itemize}
\item[(a)] the output value $J_{n,\allowbreak \nu_n,\allowbreak m_\epsilon}$ defined in
\eqref{epsilon-value-orignal} satisfies
$
|J_{n,\nu_n,m_\epsilon}-J_\nu^*|\le\epsilon;
$
\item[(b)] the output strategy $(\pi_{n,\allowbreak m_\epsilon}^1,\allowbreak \tau_{n,\allowbreak m_\epsilon}^1;
\pi_{n,\allowbreak m_\epsilon}^2,\allowbreak \tau_{n,\allowbreak m_\epsilon}^2)$ defined in
\eqref{epsilon-Nash-orignal} is an $\epsilon$-Nash equilibrium of the
 game $\mathcal G$;
\item[(c)] let $\{\epsilon_k\}_{k\ge1}$ satisfy
$\epsilon_k\downarrow0$,
$n_k\ge\lceil24\nu(w)/\epsilon_k\rceil$.
Let $(\pi_{n_k,\allowbreak m_k}^1,\allowbreak \tau_{n_k,\allowbreak m_k}^1;
\pi_{n_k,\allowbreak m_k}^2,\allowbreak \tau_{n_k,\allowbreak m_k}^2)$ be the output of
Algorithm~\ref{alg:truncVI} for $(n_k,\allowbreak \epsilon_k)$.
Then every accumulation point of
$\{(\pi_{n_k,\allowbreak m_k}^1,\allowbreak \tau_{n_k,\allowbreak m_k}^1;
\pi_{n_k,\allowbreak m_k}^2,\allowbreak \tau_{n_k,\allowbreak m_k}^2)\}_{k\ge1}$ is a Nash equilibrium of
$\mathcal G$ with respect to $(\nu,\allowbreak \lambda^*)$.
\end{itemize}
\end{theorem}

 Before proving  Theorem~\ref{thm:main_epsilon}, we need some preliminaries. Precisely, for the convergence in Theorem~\ref{thm:main_epsilon} (c) we need to  introduce suitable  topologies on the spaces $\Pi^k\times\Gamma (k=1,2)$ and then establish their compactness.
Since the policy spaces $\Pi^1$ and $\Pi^2$ are compact
and metrizable in the product topology, it remains to
give a suitable topology for the stopping-time space $\Gamma$ and then prove the compactness of $\Gamma$. Although 
the compactness of stopping times of some adapted continuous-time stochastic processes is proved \citep{BaxterChacon1977,Meyer1978}, 
we are not aware of a discrete-time version.
Inspired by the one-to-one correspondence between stopping times and
$\{0,1\}$-valued functions in \citet{Fischer2013Stopping}, we also establish a topology and the compactness of $\Gamma$, the details of which are given below. 

Let  $\mathcal H:=\bigcup_{t\ge0}H_t$ and consider functions
$d:\mathcal H\to\{0,1\}$, where $d(h_t)=1$ means that
stopping has occurred by time $t$, and $d(h_t)=0$ means
that stopping has not yet occurred.
These functions  satisfy the following: for every $t\ge0$ and $h_t\in H_t$,
if $d(h_t)=1$, then $d(h_{t+1})=1$ for every extension
$h_{t+1}=(h_t,a_t,b_t,i_{t+1},\lambda_{t+1})\in H_{t+1}$.
Denote the space of all such functions $d$ by $\mathbb D$.
Since $S$ and $\Lambda^*$ are
countable,  $A(i)$ and $B(i)(i\in S)$  are finite, $\mathcal H$ is countable, and then $\{0,1\}^{\mathcal H}$ is compact and metrizable.

To build a topology for $\Gamma$, we give two maps $\Theta$ and $\Upsilon$ as follows. $\Theta:\Gamma\to\mathbb D$ and
$\Upsilon:\mathbb D\to\Gamma$, where
\begin{align*}
\Theta(\tau)(h_t)
&:=\mathbf 1_{z_t(\{\tau\le t\})}(h_t), \ h_t\in H_t \ {\rm and} \ t\geqq 0 \quad(\text{where} \ z_t(C):=\{z_t(\omega):\omega\in C\});\\
\Upsilon(d)(\omega)
&:=\inf\{t\ge0:d(z_t(\omega))=1\},
\quad \omega\in\Omega.
\end{align*}
We equip $\mathbb D$ with the subspace topology inherited
from the product space $\{0,1\}^{\mathcal H}$, and equip
$\Gamma$ with the topology induced by $\Theta$, that is, 
$
\tau_m\to\tau \text{ in }\Gamma $ (under this topology) {\it if and only if} $\Theta(\tau_m)\to\Theta(\tau) \text{ in }\mathbb D.
$

\begin{lemma}
\label{lem:strategy-space-compact}
The map $\Theta:\Gamma\to\mathbb D$ is a homeomorphism
with inverse $\Upsilon$.
Moreover, $\Gamma$ and
$(\Pi^1\times\Gamma)\times(\Pi^2\times\Gamma)$
are compact and metrizable.
\end{lemma}
\begin{proof}
	First,  we show that $\mathbb D$ is closed.
Let $\{d_m\}_{m\ge1}\subset\mathbb D$ converge to
$d\in\{0,1\}^{\mathcal H}$, i.e.,
$
d_m(h)\longrightarrow d(h)$ for all $h\in\mathcal H.$
Suppose that $d(h_t)=1$ for some $t\ge0$ and $h_t\in H_t$.
Since each coordinate takes values in $\{0,1\}$,
$d_m(h_t)=1$ for all sufficiently large $m$.
For every extension
$h_{t+1}=(h_t,a_t,b_t,i_{t+1},\lambda_{t+1})\in H_{t+1}$,
the defining property of $\mathbb D$ gives
$d_m(h_{t+1})=1$ for all such $m$.
Taking limits yields $d(h_{t+1})=1$.
Hence $d\in\mathbb D$, and $\mathbb D$ is closed in $\{0,\allowbreak 1\}^{\mathcal H}$.
Since
$\{0,\allowbreak 1\}^{\mathcal H}$ is compact and metrizable, so is $\mathbb D$.

Next, we verify that $\Theta$ and $\Upsilon$
are mutually inverse.
If $\Theta(\tau)(h_t)=1$, then, since
$\{\tau\le t\}\in\sigma(z_t)$,
$\tau(\omega)\le t$ for every $\omega$
with $z_t(\omega)=h_t$.
Hence $\Theta(\tau)(h_{t+1})=1$ for every extension
$h_{t+1}\in H_{t+1}$ of $h_t$.
Thus $\Theta(\Gamma)\subseteq\mathbb D$.
Conversely, for $d\in\mathbb D$,
$
\{\Upsilon(d)\le t\}
=\bigcup_{s=0}^t\{d(z_s)=1\}
\in\sigma(z_t),$
so $\Upsilon(d)\in\Gamma$.
Moreover,
$
\Upsilon(\Theta(\tau))(\omega)
=\inf\{t\ge0:\tau(\omega)\le t\}
=\tau(\omega).
$
For $h_t\in H_t$, let $h_s$ denote its prefix
up to time $s\le t$.
By the defining property of $\mathbb D$,
$
\Theta(\Upsilon(d))(h_t)
=\mathbf 1_{\{\exists\,s\le t:\,d(h_s)=1\}}
=d(h_t).
$
Thus $\Theta$ and $\Upsilon$ are mutually inverse.
By the definition of the topology on $\Gamma$,
the bijection $\Theta$ is a homeomorphism.
Since $\mathbb D$ is compact and metrizable,
so is $\Gamma$.

For each $k=1,2$, the policy space $\Pi^k$ is a countable
product of compact finite-dimensional simplices.
Therefore, the full strategy space
$(\Pi^1\times\Gamma)\times(\Pi^2\times\Gamma)$
is also compact and metrizable.
\end{proof}

Note that, by the definition of the topology in $\Gamma$ above, we see that $\tau_m\to\tau$ in $\Gamma$ if and only if
\begin{equation}\label{eq:stop-convergence}
\mathbf 1_{z_t(\{\tau_m\le t\})}(h_t)
\longrightarrow
\mathbf 1_{z_t(\{\tau\le t\})}(h_t) \  \  {\rm for \ every } \ h_t\in H_t \ {\rm and} \  t\ge0.
\end{equation}


Our proof of the continuity of $J_\nu$ is to approximate it by the following finite-horizon case. Given any $T>0$, for $(\pi^1,\allowbreak \tau^1;\pi^2,\allowbreak \tau^2)\in(\Pi^1\times\Gamma)\times(\Pi^2\times\Gamma)$ and $
\tau:=\allowbreak \tau^1\wedge\tau^2\wedge\tau_D,
$ we define
\begin{eqnarray}\label{eq:finite-horizon-criterion}
\mspace{-21mu}J_{\nu,T}(\pi^1,\tau^1;\pi^2,\tau^2)
:=
\mathbb P_{\nu\otimes\delta_{\lambda^*}}^{\pi^1,\pi^2}
\bigl(X(\tau^1,\tau^2)\ge\lambda^*,\;\tau\le T\bigr),
\end{eqnarray}

\begin{lemma}\label{lem:finite-horizon-continuous}
\begin{itemize}
\item[(a)] For every fixed $T\in\mathbb N$, the functional
$J_{\nu,\allowbreak T}$ is continuous on
$(\Pi^1\times\Gamma)\times(\Pi^2\times\Gamma)$.
\item[(b)] Under the Lyapunov condition \eqref{eq:drift_section} with $
	\nu(w)<\infty$, $J_\nu
$ is continuous on $(\Pi^1\times\Gamma)\times(\Pi^2\times\Gamma)$.
\end{itemize}
\end{lemma}
\begin{proof}
(a) Let $\{(\pi_m^1,\allowbreak \tau_m^1;\pi_m^2,\allowbreak \tau_m^2)\}_{m\ge1}$ be a sequence in $(\Pi^1\times\Gamma)\times(\Pi^2\times\Gamma)$ and
converges to $(\pi^1,\allowbreak \tau^1;\pi^2,\allowbreak \tau^2)$.

For each $h_T\in H_T$, let
$$
p(h_T)
:=
\mathbb P_{\nu\otimes\delta_{\lambda^*}}^{\pi^1,\pi^2}
(z_T=h_T),~~~p_m(h_T)
:=
\mathbb P_{\nu\otimes\delta_{\lambda^*}}^{\pi_m^1,\pi_m^2}
(z_T=h_T),$$
\[
\begin{aligned}
e(h_T)&:=\mathbf 1_{z_T(\{X(\tau^1,\tau^2)\ge\lambda^*,\;\tau\le T\})}(h_T),~~
e_m(h_T)&:=\mathbf 1_{z_T(\{X(\tau_m^1,\tau_m^2)\ge\lambda^*,\tau_m\le T\})}(h_T),
\end{aligned}
\]
with $\tau:=\allowbreak \tau^1\wedge\tau^2\wedge\tau_D$ and $\tau_m:=\allowbreak \tau_m^1\wedge\tau_m^2\wedge\tau_D$.
Thus, we have
\begin{eqnarray}\label{eq:JnuT-history-sum}
J_{\nu,T}(\pi^1,\tau^1;\pi^2,\tau^2)
=
\sum_{h_T\in H_T}
e(h_T)\,p(h_T),
\end{eqnarray}
\begin{eqnarray}\label{eq:JnuTk-history-sum}
J_{\nu,T}(\pi_m^1,\tau_m^1;\pi_m^2,\tau_m^2)
=
\sum_{h_T\in H_T}
e_m(h_T)\,p_m(h_T).
\end{eqnarray}
Let $\pi^k_m=\allowbreak :(\pi_{m,\allowbreak t}^k,\allowbreak  t\geq 0) (k=\allowbreak 1,\allowbreak 2)$ for $m\geq 1$. Then,  $\lim_{m\to\infty}\pi_{m,\allowbreak t}^k(\cdot|h_t)=\allowbreak \pi_t^k(\cdot|h_t)$ for   $h_t\in H_t$ and $t\ge0$, thus
\begin{eqnarray}\label{eq:total-variation-convergence}
 \lim_{m\to\infty}p_m(h_T)=p(h_T) \ \ \forall h_T\in H_T,  \ {\rm as} \ m\to\infty.
\end{eqnarray}

Fix any $h_T=\allowbreak (i_0,\allowbreak \lambda_0,\allowbreak a_0,\allowbreak b_0,\allowbreak \dots,\allowbreak i_T,\allowbreak \lambda_T)\in H_T$, and for each $t=\allowbreak 0,\allowbreak 1,\allowbreak \dots,\allowbreak T$ let
$
h_t:=\allowbreak (i_0,\allowbreak \lambda_0,\allowbreak a_0,\allowbreak b_0,\allowbreak \dots,\allowbreak i_t,\allowbreak \lambda_t).
$
Then $h_t\in H_t$ and $h_{t+1}=\allowbreak (h_t,\allowbreak a_t,\allowbreak b_t,\allowbreak i_{t+1},\allowbreak \lambda_{t+1})$ for every $t<T$.
By \eqref{eq:stop-convergence},
$$
\mathbf 1_{z_t(\{\tau_m^k\le t\})}(h_t)
\longrightarrow
\mathbf 1_{z_t(\{\tau^k\le t\})}(h_t)
~~ \forall\, t=0,1,\dots,T,\ k=1,2.
$$
Since these indicators take only values $0$ or $1$, there exists $M\in\mathbb N$ such that
$$
\mathbf 1_{z_t(\{\tau_m^k\le t\})}(h_t)
=
\mathbf 1_{z_t(\{\tau^k\le t\})}(h_t)
~~\forall\,m\ge M,\ t=0,\dots,T,\ k=1,2.
$$

By the convention $\inf\emptyset=\allowbreak \infty$, the equality of the above
indicators implies that, for every $m\ge M$ and $k=\allowbreak 1,\allowbreak 2$,
\begin{eqnarray}\label{eq:tau-from-indicators}
\inf\left\{t\in\{0,\ldots,T\}:
\mathbf 1_{z_t(\{\tau_m^k\le t\})}(h_t)=1\right\}
=
\inf\left\{t\in\{0,\ldots,T\}:
\mathbf 1_{z_t(\{\tau^k\le t\})}(h_t)=1\right\}.
\end{eqnarray}
It follows that, for all $m\ge M$,
$
\mathbf 1_{z_T(\{\tau_m\le T\})}(h_T)
=\allowbreak
\mathbf 1_{z_T(\{\tau\le T\})}(h_T)$, where $\tau_m:=\allowbreak \tau_m^1\wedge\tau_m^2\wedge\tau_D$ and $\tau:=\allowbreak \tau^1\wedge\tau^2\wedge\tau_D$.\\
If
$
\mathbf 1_{z_T(\{\tau_m\le T\})}(h_T)
=\allowbreak
\mathbf 1_{z_T(\{\tau\le T\})}(h_T)=\allowbreak 1$, then, for any $\omega\in\Omega$ satisfying $z_T(\omega)=h_T$, we have
$\tau_m(\omega)=\tau(\omega)$ and
$X(\tau_m^1,\tau_m^2)(\omega)=X(\tau^1,\tau^2)(\omega)$.
By the definitions of $e_m(h_T)$ and $e(h_T)$ using the joint events above, it follows that $e_m(h_T)=e(h_T)$.\\
If $
\mathbf 1_{z_T(\{\tau_m\le T\})}(h_T)=\allowbreak
\mathbf 1_{z_T(\{\tau\le T\})}(h_T)=\allowbreak 0$, then $e_m(h_T)=\allowbreak e(h_T)=\allowbreak 0$.
Thus, $e_m(h_T)=\allowbreak e(h_T)$ for all $m\ge M$.
Hence
$$\lim_{m\to\infty} e_m(h_T)=e(h_T) \ \ \forall h_T\in H_T,$$
which, together with (\ref{eq:JnuT-history-sum})-(\ref{eq:total-variation-convergence}) and  Lemma 8.3.7 in \citet{Hernandez1999Further}, implies that
$\lim_{m\to\infty} J_{\nu,\allowbreak T}(\pi_m^1,\allowbreak \tau_m^1;\pi_m^2,\allowbreak \tau_m^2)=\allowbreak J_{\nu,\allowbreak T}(\pi^1,\allowbreak \tau^1;\pi^2,\allowbreak \tau^2)$, proving the continuity.\\
(b)
By the definition of $J_\nu$ and $J_{\nu,\allowbreak T}$, under the Lyapunov condition
\eqref{eq:drift_section}, we have
\begin{eqnarray*}
|
J_\nu(\pi^1,\tau^1;\pi^2,\tau^2)
-
J_{\nu,T}(\pi^1,\tau^1;\pi^2,\tau^2)
|
\le
\mathbb P_{\nu\otimes\delta_{\lambda^*}}^{\pi^1,\pi^2}(\tau>T)
\end{eqnarray*}
 for $(\pi^1,\allowbreak \tau^1;\pi^2,\allowbreak \tau^2)\in (\Pi^1\times\Gamma)\times(\Pi^2\times\Gamma)$, with  $\tau=\allowbreak \tau^1\wedge\tau^2\wedge\tau_D$.

Moreover, by the Lyapunov condition
\eqref{eq:drift_section} and \eqref{prop:drift-2}, there exists $\rho\in(0,\allowbreak 1)$ such that
\begin{eqnarray}\label{eq:lyapunov-tail}
\mathbb P_{(i,\lambda)}^{\pi^1,\pi^2}(\tau>T)
\le\mathbb P_{(i,\lambda)}^{\pi^1,\pi^2}(\tau_D>T)
\le
\rho^T w(i),
\end{eqnarray}
for   $(i,\allowbreak \lambda)\in D^c\times\Lambda^*$ and $(\pi^1,\allowbreak \tau^1;\pi^2,\allowbreak \tau^2)\in(\Pi^1\times\Gamma)\times(\Pi^2\times\Gamma)$.

Integrating \eqref{eq:lyapunov-tail} with respect to $\nu$ gives
\begin{eqnarray}\label{eq:uniform-tail-bound}
\mspace{-21mu}\sup_{(\pi^1,\tau^1;\pi^2,\tau^2)}
\mathbb P_{\nu\otimes\delta_{\lambda^*}}^{\pi^1,\pi^2}(\tau>T)
\le
\rho^T\nu(w) \ \ \ \forall T>0.
\end{eqnarray}
Now, let
$(\pi_m^1,\allowbreak \tau_m^1;\pi_m^2,\allowbreak \tau_m^2)\to(\pi^1,\allowbreak \tau^1;\pi^2,\allowbreak \tau^2)$.   Then, by (\ref{eq:uniform-tail-bound}) we have
\begin{eqnarray*}
&&\bigl|
J_\nu(\pi_m^1,\tau_m^1;\pi_m^2,\tau_m^2)
-
J_\nu(\pi^1,\tau^1;\pi^2,\tau^2)
\bigr|
\\
&\le&
\bigl|
J_\nu(\pi_m^1,\tau_m^1;\pi_m^2,\tau_m^2)
-
J_{\nu,T}(\pi_m^1,\tau_m^1;\pi_m^2,\tau_m^2)
\bigr|
+
\bigl|
J_{\nu,T}(\pi_m^1,\tau_m^1;\pi_m^2,\tau_m^2)
-
J_{\nu,T}(\pi^1,\tau^1;\pi^2,\tau^2)
\bigr|\\
&&
+
\bigl|
J_{\nu,T}(\pi^1,\tau^1;\pi^2,\tau^2)
-
J_\nu(\pi^1,\tau^1;\pi^2,\tau^2)
\bigr|\\
&
\leq& 2\rho^T\nu(w)+\bigl|
J_{\nu,T}(\pi_m^1,\tau_m^1;\pi_m^2,\tau_m^2)
-
J_{\nu,T}(\pi^1,\tau^1;\pi^2,\tau^2)
\bigr|,
\end{eqnarray*}
which, together with the continuity of $J_{\nu,\allowbreak T}$, gives
$\lim_{m\to\infty}J_\nu(\pi_m^1,\allowbreak \tau_m^1;\pi_m^2,\allowbreak \tau_m^2)=\allowbreak J_\nu(\pi^1,\allowbreak \tau^1;\pi^2,\allowbreak \tau^2).$
\end{proof}

\begin{proof}[Proof of Theorem~\ref{thm:main_epsilon}]
(a) By Theorem~\ref{thm:truncation}(a) and $n\ge \frac{24\nu(w)}{\epsilon}$, we have
\begin{eqnarray}
|J_{n,\nu_n}^*-J_\nu^*|\le\frac{9\nu(w)}{n}\leq 3\epsilon/8.\label{WGGZ-2}
\end{eqnarray}
 In Step~3 of Algorithm~\ref{alg:truncVI}, the stopping rule
$
\|u_{n,\allowbreak m_\epsilon}^1-u_{n,\allowbreak m_\epsilon}^0\|_{E_n}\le\epsilon/4
$
and $u_{n,\allowbreak m_\epsilon}^0\le u_n^*\le u_{n,\allowbreak m_\epsilon}^1$ imply
\[
|u_{n,m_\epsilon}^0(i,\lambda)-u_n^*(i,\lambda)|\le\frac{\epsilon}{4}
\quad \forall (i,\lambda)\in E_n.
\]
Therefore,
\[
|J_{n,\nu_n,m_\epsilon}-J_{n,\nu_n}^*|
=
\Bigl|
\sum_{i\in D_n^c}u_{n,m_\epsilon}^0(i,\lambda^*)\nu_n(i)
-
J_{n,\nu_n}^*
\Bigr|
\le\frac{\epsilon}{4},
\]
which,together with \eqref{WGGZ-2}, gives
\[
|J_{n,\nu_n,m_\epsilon}-J_\nu^*|
\le
|J_{n,\nu_n,m_\epsilon}-J_{n,\nu_n}^*|
+
|J_{n,\nu_n}^*-J_\nu^*|
\leq\epsilon.
\]
Hence, (a) is true.

(b) By Theorem~\ref{thm:truncated-convergence}(a), we see that
$\widehat{s}_{n,\allowbreak m_\epsilon}:=\allowbreak
\bigl(
\widehat{\pi}_{n,\allowbreak m_\epsilon}^1,\allowbreak \widehat{\tau}_{n,\allowbreak m_\epsilon}^1;
\widehat{\pi}_{n,\allowbreak m_\epsilon}^2,\allowbreak \widehat{\tau}_{n,\allowbreak m_\epsilon}^2
\bigr)$ is an $\epsilon/4$-Nash equilibrium of
$\mathcal G_n$ with respect to $(\nu_n,\allowbreak \lambda^*)$.
By Theorem~\ref{thm:truncation}(b), $\Psi(\widehat{s}_{n,\allowbreak m_\epsilon})$ is an
$(\epsilon/4+18\nu(w)/n)$-Nash equilibrium of the original game with
respect to $(\nu,\allowbreak \lambda^*)$.
Since $n\ge24\nu(w)/\epsilon$,
\[
\frac{\epsilon}{4}+\frac{18\nu(w)}{n}
\le
\frac{\epsilon}{4}+\frac{18\epsilon}{24}
=
\frac{\epsilon}{4}+\frac{3\epsilon}{4}
=
\epsilon,
\]
$\Psi(\widehat{s}_{n,\allowbreak m_\epsilon})$ is an $\epsilon$-Nash equilibrium of
$\mathcal G$, that is, (b) is proved.

(c) By (a)-(b) just proved,$(\pi_{n_k,\allowbreak m_k}^1,\allowbreak \tau_{n_k,\allowbreak m_k}^1;\pi_{n_k,\allowbreak m_k}^2,\allowbreak \tau_{n_k,\allowbreak m_k}^2)$ is an $\epsilon_k$-Nash equilibrium of $\mathcal G.$ Thus, by Definition~\ref{def:epsilon}, we see that, for each $(\pi^1,\allowbreak \tau^1)\in\Pi^1\times\Gamma$ and $(\pi^2,\allowbreak \tau^2)\in\Pi^2\times\Gamma$,
\begin{eqnarray}
J_\nu(\pi_{n_k,m_k}^1,\tau_{n_k,m_k}^1;\pi^2,\tau^2)
&\ge &J_\nu^*-\epsilon_k,
\label{eq:eps-opt-1}\\
J_\nu(\pi^1,\tau^1;\pi_{n_k,m_k}^2,\tau_{n_k,m_k}^2)
&\le& J_\nu^*+\epsilon_k.
\label{eq:eps-opt-2}
\end{eqnarray}
On the other hand,  Lemma~\ref{lem:strategy-space-compact} guarantees the compactness of $(\Pi^1\times\Gamma)\times(\Pi^2\times\Gamma)$. Hence, the sequence
$\{(\pi_{n_k,\allowbreak m_k}^1,\allowbreak \tau_{n_k,\allowbreak m_k}^1;
\pi_{n_k,\allowbreak m_k}^2,\allowbreak \tau_{n_k,\allowbreak m_k}^2)\}_{k\ge1}$
has an accumulation point, denoted by
$(\pi^1_*,\allowbreak \tau^1_*;\pi^2_*,\allowbreak \tau^2_*)$.
Passing to a subsequence and keeping the same index notation, we may assume
\begin{eqnarray}
\lim_{k\to\infty}(\pi_{n_k,m_k}^1,\tau_{n_k,m_k}^1;
\pi_{n_k,m_k}^2,\tau_{n_k,m_k}^2)=(\pi^1_*,\tau^1_*;\pi^2_*,\tau^2_*).\label{WGGZ-3}
\end{eqnarray}
Moreover, fix  $(\pi^1,\allowbreak \tau^1)\in\Pi^1\times\Gamma$. Taking $\pi^2=\allowbreak \pi_{n_k,\allowbreak m_k}^2$ and $\tau^2=\allowbreak \tau_{n_k,\allowbreak m_k}^2$ in \eqref{eq:eps-opt-1}, and using \eqref{eq:eps-opt-2}, we obtain
\begin{eqnarray}
J_\nu(\pi_{n_k,m_k}^1,\tau_{n_k,m_k}^1;\pi_{n_k,m_k}^2,\tau_{n_k,m_k}^2)\ge J_\nu^*-\epsilon_k
\ge J_\nu(\pi^1,\tau^1;\pi_{n_k,m_k}^2,\tau_{n_k,m_k}^2)-2\epsilon_k.\label{WGGZ-4}
\end{eqnarray}
Thus, using Lemma~\ref{lem:finite-horizon-continuous} and $\lim_{k\to\infty}\epsilon_k=\allowbreak 0$, by (\ref{WGGZ-4}) we have
\[
J_\nu(\pi^1_*,\tau^1_*;\pi^2_*,\tau^2_*)
\ge
J_\nu(\pi^1,\tau^1;\pi^2_*,\tau^2_*).
\]
Similarly, fix  $(\pi^2,\allowbreak \tau^2)\in\Pi^2\times\Gamma$.
Taking $\pi^1=\allowbreak \pi_{n_k,\allowbreak m_k}^1$ and $\tau^1=\allowbreak \tau_{n_k,\allowbreak m_k}^1$ in \eqref{eq:eps-opt-2}, and using \eqref{eq:eps-opt-1}, we obtain
\begin{eqnarray*}
J_\nu(\pi_{n_k,m_k}^1,\tau_{n_k,m_k}^1;
\pi_{n_k,m_k}^2,\tau_{n_k,m_k}^2)
\le J_\nu^*+\epsilon_k
\le J_\nu(\pi_{n_k,m_k}^1,\tau_{n_k,m_k}^1;\pi^2,\tau^2)+2\epsilon_k.
\end{eqnarray*}
Letting $k\to\infty$ gives
$
J_\nu(\pi^1_*,\tau^1_*;\pi^2_*,\tau^2_*)
\le
J_\nu(\pi^1_*,\tau^1_*;\pi^2,\tau^2).
$
Thus $(\pi^1_*,\allowbreak \tau^1_*;\pi^2_*,\allowbreak \tau^2_*)$ is a Nash equilibrium of $\mathcal G$
with respect to $(\nu,\allowbreak \lambda^*)$.
\end{proof}
\begin{remark}\label{rem:algorithm-convergence-novelty}
Under a weaker Lyapunov condition,
Theorem~\ref{thm:truncated-convergence} shows that, for each fixed
truncation level and any prescribed accuracy, Algorithm~\ref{alg:truncVI} obtains an approximate
value and an approximate Nash equilibrium of the truncated game in finitely many iterations.
Theorem~\ref{thm:main_epsilon} establishes the corresponding approximate results
for the original game with countably infinite states.
By contrast, the algorithms in
\citet{HGW-2023,GuoWen2025Nonstationary,Wen2022Optimal}
generally require finite-state representations for numerical
implementation; those in
\citet{GuoWen2025Nonstationary,Wen2022Optimal}
also rely on conditions of the form \eqref{one-step-absorb}. Besides, we  further provide convergence analysis of the approximate values and  Nash equilibria for the truncated  and original games, which are new.
\end{remark}

\section{An Example}\label{sec:numerics}

In this section, we apply our main results to a controlled queueing
system, and illustrate the implementation and numerical behavior of the
proposed algorithm. We will see that the controlled queueing system satisfies the Lyapunov condition \eqref{eq:drift_section} but
violates the condition \eqref{one-step-absorb} used commonly in the existing literature.

\subsection*{Example 7.1 (A controlled queueing system)}

Consider a single-server queueing system in discrete time, with a state space
$S=\allowbreak \{0,\allowbreak 1,\allowbreak 2,\allowbreak \ldots\}$ representing the set of numbers of customers at the beginning of each slot. A server of this system chooses a service capacity, determining the service-completion probability; an arrival controller chooses an arrival level, determining the arrival probability.
 The system  is naturally closed once there is not any customer (i.e., the state is at 0).

We view the server as Player 1 and the arrival controller as Player 2. At each decision epoch, either player may stop the game early, or both may continue.
Suppose that the state (i.e., the queue length) at decision epoch $t$ is $i_t\ge1$.
If Player 1 quits the system, the terminal payoff is $g^1(i_t)$;
if Player 2 quits  while Player 1 continues, the terminal payoff is $g^2(i_t)$;
if both continue, Player 1 chooses a service mode
    $a_t$ from a set $A(i_t)$ consisting of possible service modes at state $i_t$ for player 1, and player~2 chooses an arrival mode
    $b_t$ from a set $B(i_t)$ of available modes, and Player 1 receives the one-stage payoff $r(i_t,\allowbreak a_t,\allowbreak b_t)$ from player~2.

 Since their interests are directly opposed, this is modeled as a two-person zero-sum stochastic game. Player 1 aims to maximize the probability that the accumulated payoff reaches a threshold
  $\lambda^*>0$, while Player 2 aims to minimize it.

In this example, suppose that the service and arrival modes consist of low $(
L)$ or high $(H)$ modes. Thus, the game is formally specified as follows:
\begin{enumerate}[(i)]
\item The state space is $S=\allowbreak \{0,\allowbreak 1,\allowbreak \ldots,\allowbreak i,\allowbreak \ldots\}$, and the target set is $D=\allowbreak \{0\}$.
\item For each $i\in S$, $A(i)=\allowbreak B(i)=\allowbreak \{L,\allowbreak H\},\allowbreak $ where $L$ and $H$ denote the low and high modes, respectively.
\item For each $i\ge1$, let $\mu_i(a)>0$ denote the probability of one
    service completion during a time slot under service action $a$, and
    let $\lambda_i(b)\ge0$ denote the probability of one arrival
    during that slot under an arrival action $b$. Assume that
    $
    \mu_i(a)+\lambda_i(b)\le1,\allowbreak
    ~~
    i\ge1,\allowbreak ~ a\in A(i),\allowbreak ~ b\in B(i).
    $
  The transition probabilities of the game are given as follows: for each $i\in S,\allowbreak  a\in A(i),\allowbreak  b\in B(i),\allowbreak $
    \begin{eqnarray}
     p(i-1|i,a,b)=\mu_i(a),\   p(i+1|i,a,b)=\lambda_i(b), \
    p(i|i,a,b)=1-\mu_i(a)-\lambda_i(b)~~\forall~i\ge 1. \label{WGGZ-5}
    \end{eqnarray}
    Once the process enters $D$, the game terminates; equivalently,
    $p(0|0,\allowbreak a,\allowbreak b)=\allowbreak 1$ for all $a\in A(0)$ and $b\in B(0)$.

	    \item The one-stage payoff $r(i,a,b)>0$ and the terminal payoffs $g^1(i)>0$ and
	    $g^2(i)>0$ for all $i\ge1$ and $a\in A(i)$ and $b\in B(i)$.
\end{enumerate}

To verify the Lyapunov condition in
Proposition~\ref{prop:attenuation}, we impose the following condition.
\begin{condition}\label{cond:queue-stability}
The service  and arrival probabilities satisfy
\begin{eqnarray}\label{probability-condition}
\sup_{\substack{i\ge1,\,b\in B(i)}}\lambda_i(b)<\inf_{\substack{i\ge1,\,a\in A(i)}}\mu_i(a).
\end{eqnarray}
\end{condition}

\begin{proposition}
Under Condition~\ref{cond:queue-stability}, the controlled queueing
game  in Example~7.1 satisfies the Lyapunov condition
\eqref{eq:drift_section}. Consequently, the game has a value and a
Nash equilibrium.
\end{proposition}

\begin{proof}
We verify the Lyapunov condition in
Proposition~\ref{prop:attenuation}.
If
$
\sup_{\substack{i\ge1,\allowbreak \,b\in B(i)}}\lambda_i(b)=\allowbreak 0,\allowbreak
$
choose any $\theta>1$.
If
$
\sup_{\substack{i\ge1,\allowbreak \,b\in B(i)}}\lambda_i(b)>0,\allowbreak
$
choose
\begin{eqnarray}\label{theta-chose}
1<\theta<
\frac{\inf_{i\ge1,\,a\in A(i)}\mu_i(a)}
{\sup_{i\ge1,\,b\in B(i)}\lambda_i(b)}.
\end{eqnarray}
In either case, let
$
w(i)=\allowbreak \theta^i.
$
For $i\ge2$, by (\ref{WGGZ-5}) and (\ref{theta-chose}) we have
\begin{eqnarray*}
&&\sum_{j\in D^c}p(j|i,a,b)w(j)
=\mu_i(a)\theta^{i-1}
+\lambda_i(b)\theta^{i+1}
+\bigl(1-\mu_i(a)-\lambda_i(b)\bigr)\theta^i\\
&=&\theta^i
\left[
1+(\theta-1)
\left(
\lambda_i(b)-\frac{\mu_i(a)}{\theta}
\right)
\right]\le\theta^i
\left[
1+(\theta-1)\big(\sup_{k,d}\lambda_k(d)
-\frac{1}{\theta}
\inf_{k,c}\mu_k(c)\big)
\right].
\end{eqnarray*}
By the choice of $\theta$,
$
\theta
\sup_{k,\allowbreak d}\lambda_k(d)
<
\inf_{k,\allowbreak c}\mu_k(c),\allowbreak
$
where the inequality is immediate if $\sup_{k,\allowbreak d}\lambda_k(d)=\allowbreak 0$. Hence
\[
\rho_1:=1+(\theta-1)
\left[
\sup_{k,d}\lambda_k(d)
-\frac{1}{\theta}
\inf_{k,c}\mu_k(c)
\right]
<1.
\]
Also, for $i=\allowbreak 1$, using (\ref{WGGZ-5}) and (\ref{theta-chose}) again, we have
\begin{eqnarray*}
&&\sum_{j\in D^c}p(j\mid1,a,b)w(j)
=\bigl(1-\mu_1(a)-\lambda_1(b)\bigr)\theta
+\lambda_1(b)\theta^2\\
&=&
\theta
\left[
1-\mu_1(a)+\lambda_1(b)(\theta-1)
\right]
\le
\theta
\left[
1-
\inf_{k,c}\mu_k(c)
+(\theta-1)
\sup_{k,d}\lambda_k(d)
\right].
\end{eqnarray*}
By the choice of $\theta$, we have 
$
(\theta-1)
\sup_{k,\allowbreak d}\lambda_k(d)
<
\inf_{k,\allowbreak c}\mu_k(c)\allowbreak.$ Therefore,
$
\rho_2:=1-
\inf_{k,c}\mu_k(c)
+(\theta-1)
\sup_{k,d}\lambda_k(d)
<1.
$
Consequently, for the constant $\rho:=\allowbreak \rho_1\vee \rho_2<1$, we have
$
\sum_{j\in D^c}p(j|i,a,b)w(j)\le\rho w(i), \ \ i\ge1,a\in A(i),b\in B(i).
$
Thus, the Lyapunov condition \eqref{eq:drift_section} holds.
Theorem~\ref{game-Nash} then ensures the existence of the value and a
Nash equilibrium of this queueing game.
\end{proof}

\begin{remark}
Condition~\ref{cond:queue-stability} implies the Lyapunov
condition and is therefore sufficient for
Assumption~\ref{ass:absorb}. However, condition
\eqref{one-step-absorb} used in
\citet{HG-2020,GuoWen2025Nonstationary,Wen2022Optimal,Liu2001,
Huang2011} fails in this example. This is because $p(0|i,\allowbreak a,\allowbreak b)=\allowbreak 0$ for each  $(a,\allowbreak b)\in A(i)\times B(i)$ with $i\geq 2$.
\end{remark}

Now we perform numerical calculations for this queueing game.
For the numerical experiment, we take
\begin{eqnarray*}
\mu_i(L)=0.375,~~
\mu_i(H)=0.458,~~
\lambda_i(L)=0.092,~~
\lambda_i(H)=0.248 ~~\forall i\ge1.
\end{eqnarray*}
It is easy to verify that the transition probabilities are well-defined and Condition~\ref{cond:queue-stability} holds.
We take the Lyapunov function
$w(i)=\allowbreak \theta^i$ with
$\theta=\allowbreak 1.256048$, which satisfies \eqref{theta-chose}.
For these data, we have
\[
\rho
=
\max
\left\{
\begin{array}{l}
1+0.248(\theta-1)
-0.375(\theta-1)/\theta,\\[1mm]
1-0.375+0.248(\theta-1)
\end{array}
\right\}
\approx 0.987<1.
\]
so that the drift condition \eqref{eq:drift_section} holds with this $\rho$.
 For every $i\ge 1$ and $a,\allowbreak b\in\{L,\allowbreak H\}$,
let the  payoffs be
\[
r(i,a,b)=
\begin{cases}
0.45, & a=b,\\
\kappa(i), & a\ne b,
\end{cases}
~~\text{with}~~
\kappa(i)=
\begin{cases}
1.20, & 1\le i\le3,\\
1.55, & 4\le i\le6,\\
1.85, & i\ge7,
\end{cases}
\]
\[
g^1(i)=
\begin{cases}
0.05, & 1\le i\le3,\\
0.25, & 4\le i\le6,\\
0.45, & i\ge7,
\end{cases}
\qquad
g^2(i)=
\begin{cases}
0.10, & i=1,\\
2.05, & 2\le i\le6,\\
0.50, & i\ge7.
\end{cases}
\]
The initial distribution and target threshold are
$
\nu(2)=0.5,
\nu(4)=0.3,
\nu(6)=0.2,
\lambda^*=2.$

Using these data  and Algorithm~1, we compute approximations of
the game value and an $\epsilon$-Nash equilibrium. For each
$\epsilon>0$, Algorithm~1 chooses
$
n=\left\lceil\frac{24\nu(w)}{\epsilon}\right\rceil,
\nu(w)=\sum_{i\in D^c}w(i)\nu(i)\approx2.321.
$
Since $w(i)=\allowbreak \theta^i$,
$
S_n\setminus D
=\allowbreak
\left\{
1,\allowbreak \ldots,\allowbreak
\left\lfloor\frac{\log n}{\log\theta}\right\rfloor
\right\}.
$
\begin{figure}[htbp]
\centering
\begin{minipage}[c]{0.48\textwidth}
    \centering
    \captionof{table}{Truncation levels, approximate values,
    and iteration counts}
    \label{tab:two_mode_convergence}

    \resizebox{\linewidth}{!}{%
    \begin{tabular}{crrrrr}
    \hline
    $\epsilon$ & $n$ & $|S_n\setminus D|$
    & $m_\epsilon$
    & $u^0_{n,m_\epsilon}(2,2)$
    & $J_{n,\nu_n,m_\epsilon}$\\
    \hline
    0.500 & 112   & 20 & 3 & 0.397 & 0.622\\
    0.250 & 223   & 23 & 4 & 0.414 & 0.658\\
    0.100 & 558   & 27 & 5 & 0.427 & 0.691\\
    0.050 & 1115  & 30 & 5 & 0.427 & 0.691\\
    0.025 & 2229  & 33 & 5 & 0.427 & 0.691\\
    0.010 & 5571  & 37 & 5 & 0.427 & 0.691\\
    0.005 & 11141 & 40 & 5 & 0.427 & 0.691\\
    0.001 & 55702 & 47 & 5 & 0.427 & 0.691\\
    \hline
    \end{tabular}%
    }
    \par\vspace{8pt}
    {\footnotesize\raggedright
    \textit{Note.} For $\epsilon=0.001$, $47$ non-target states are retained. Since $w(i)=\theta^i$,
    the retained state count grows logarithmically
    with the truncation level $n$.\par}
\end{minipage}
\hfill
\begin{minipage}[c]{0.45\textwidth}
    \centering
    \includegraphics[width=\linewidth]
    {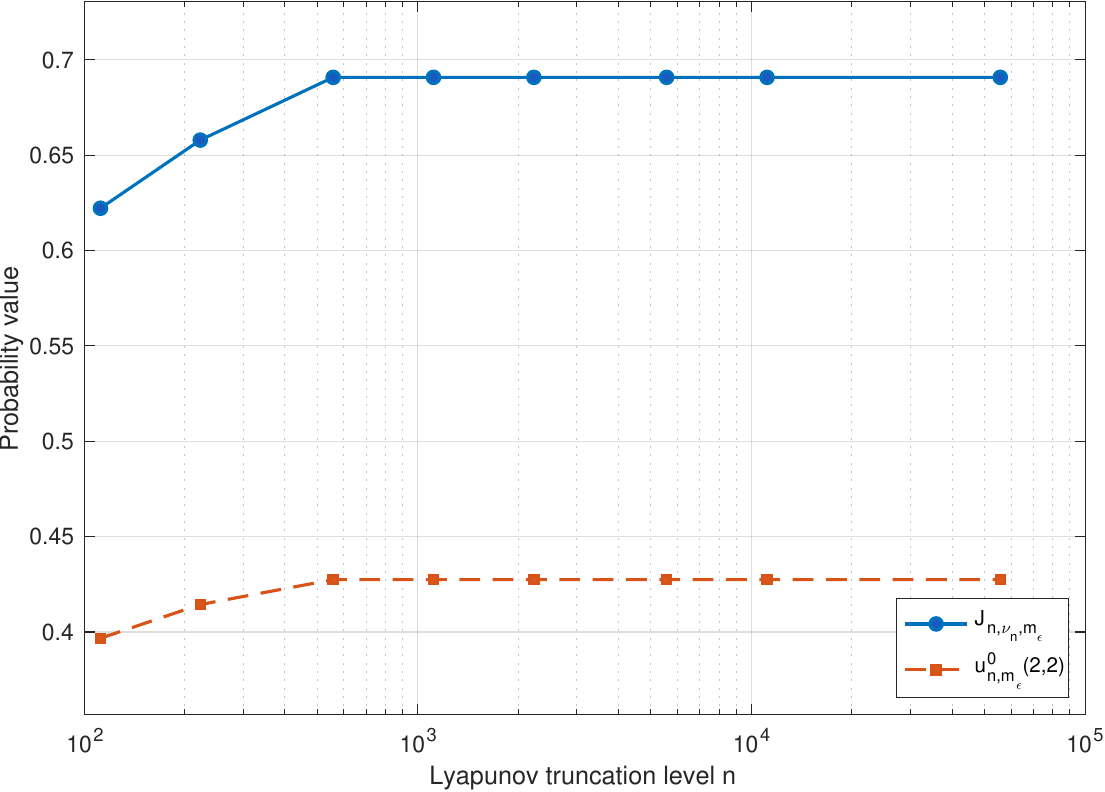}
    \captionof{figure}{The pointwise lower iterate
$u^0_{n,m_\epsilon}(2,2)$ and the distributional output
$J_{n,\nu_n,m_\epsilon}$}
    \label{fig:two_mode_value_convergence}
\end{minipage}
\end{figure}

Combining the truncation bound \eqref{eq:truncation-value-bound} with the
stopping tolerance gives
$|J_{n,\allowbreak \nu_n,\allowbreak m_\epsilon}-J_\nu^*|\le \epsilon/4+9\nu(w)/n
\le5\epsilon/8<\epsilon$. Figure~\ref{fig:two_mode_value_convergence}
shows that the two computed outputs stabilize as $\epsilon$ decreases and
the corresponding truncation level $n$ increases.

Fig.~\ref{fig:two_mode_policy_maps} displays the approximate
stationary strategies computed with $\epsilon=0.001$
for the first eight non-target states.
The upper panels show the high-mode probabilities
when the respective player continues, whereas the
lower panels show the immediate stopping regions.
Blank cells in the upper panels indicate immediate stopping.

\begin{figure}[htbp]
\centering
    \includegraphics[width=0.3\linewidth,trim=0bp 0bp 362.5bp 0bp,clip]
    {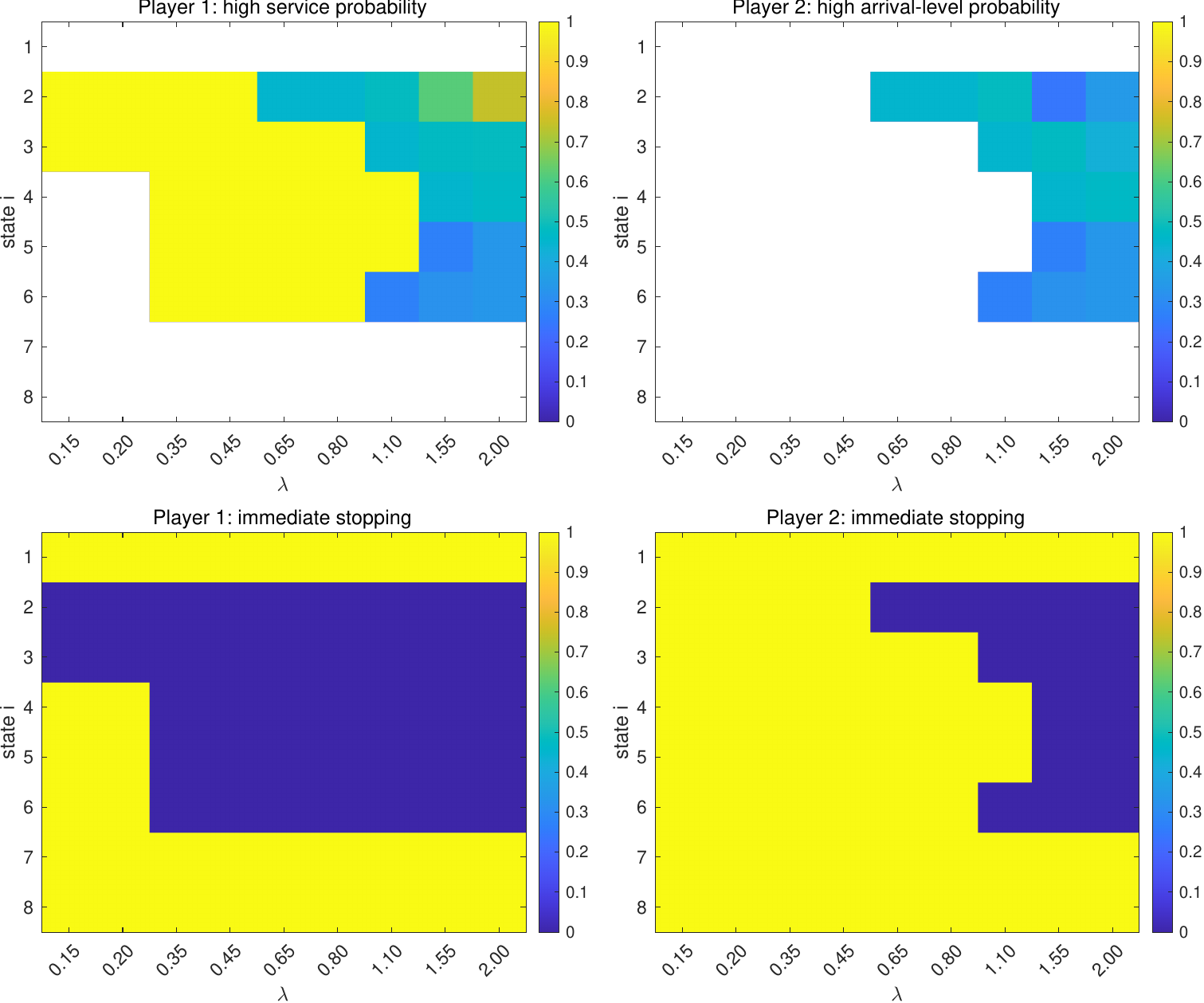}%
    \hspace{0.03\linewidth}%
    \includegraphics[width=0.3\linewidth,trim=362.5bp 0bp 0bp 0bp,clip]
    {fig_two_mode_policy_maps.pdf}
    \caption{Approximate stationary strategies for $\epsilon=0.001$.}
    \label{fig:two_mode_policy_maps}
\end{figure}

Table~\ref{tab:two_mode_strategy} reports the stopping
decisions and high-mode probabilities for selected
$(i,\lambda)$ pairs.
Here, ``Cont.'' and ``Stop'' denote continuation and
immediate stopping, respectively, and a dash indicates
that no action is selected.

\begin{table}[htbp]
    \centering
    \caption{Selected approximate stationary
    strategies for $\epsilon=0.001$}
    \label{tab:two_mode_strategy}

    \small
    \renewcommand{\arraystretch}{1.1}
    \begin{tabular}{cccccc}
    \hline
    $i$ & $\lambda$ & P1 & $\phi_H$
        & P2 & $\psi_H$\\
    \hline
    2 & 0.45 & Cont. & 1.000 & Stop  & --\\
    2 & 0.65 & Cont. & 0.450 & Cont. & 0.450\\
    2 & 2.00 & Cont. & 0.744 & Cont. & 0.345\\
    4 & 0.20 & Stop  & --    & Stop  & --\\
    4 & 1.10 & Cont. & 1.000 & Stop  & --\\
    4 & 1.55 & Cont. & 0.450 & Cont. & 0.450\\
    6 & 0.20 & Stop  & --    & Stop  & --\\
    6 & 1.10 & Cont. & 0.271 & Cont. & 0.271\\
    6 & 2.00 & Cont. & 0.337 & Cont. & 0.337\\
    \hline
    \end{tabular}
\end{table}

Overall, the example illustrates two features of the proposed method.
First, the queue satisfies the Lyapunov condition ~\eqref{eq:drift_section}, but it does not satisfy
condition \eqref{one-step-absorb} used in
\citet{HG-2020,GuoWen2025Nonstationary,Wen2022Optimal,Liu2001,
Huang2011}. Second, the exponential choice $w(i)=\allowbreak \theta^i$ reduces a
large theoretical truncation level to a comparatively small retained
state space, while the reported value approximations stabilize rapidly
at the displayed precision. Together, these observations illustrate
the effectiveness of Algorithm~1.

\setlength{\bibsep}{-0.6pt}

\end{document}